\documentclass[12pt]{article}
\usepackage{amsfonts}
\usepackage{amssymb}
\usepackage{amsmath,enumerate}
\usepackage{color}
=0pt
\def\sqr#1#2{{\vcenter{\vbox{\hrule height.#2pt
              \hbox{\vrule width.#2pt height#1pt \kern#1pt \vrule
width.#2pt}
              \hrule height.#2pt}}}}
\def\signed #1{{\unskip\nobreak\hfil\penalty50
              \hskip2em\hbox{}\nobreak\hfil#1
              \parfillskip=0pt \finalhyphendemerits=0 \par}}
\def\endpf{\signed {$\sqr69$}}
\def\dbD{{\mathbb{D}}}
\def\dbE{{\mathbb{E}}}
\def\dbF{{\mathbb{F}}}

\def\dbH{{\mathbb{H}}}

\def\dbL{{\mathbb{L}}}

\def\dbN{{\mathbb{N}}}

\def\dbP{{\mathbb{P}}}

\def\dbR{{\mathbb{R}}}
\def\dbS{{\mathbb{S}}}

\def\b{\beta}

\def\d{\delta}
\def\e{\varepsilon}

\def\t{\tau}
\def\f{\varphi}
\def\th{\theta}
\def\o{\omega}
\def\p{\phi}
\def\3n{\negthinspace \negthinspace \negthinspace }
\def\2n{\negthinspace \negthinspace }
\def\1n{\negthinspace }
\def\ns{\noalign{\smallskip} }

\def\ds{\displaystyle}
\def\O{\Omega}
\def\Om{\Omega}
\def\cA{{\cal A}}

\def\cD{{\cal D}}

\def\cF{{\cal F}}

\def\cJ{{\cal J}}

\def\cL{{\cal L}}

\def\cQ{{\cal Q}}

\def\cU{{\cal U}}

\def\mE{{\mathbb{E}}}

\def\no{\noindent}

\def\ms{\medskip}
\def\bs{\bigskip}
\def\q{\quad}
\def\qq{\qquad}
\def\hb{\hbox}

\def\liminf{\mathop{\underline{\rm lim}}}

\def\lan{\mathop{\langle}}

\def\pa{\partial}

\def\wt{\widetilde}
\def\cd{\cdot}
\def\cds{\cdots}

\def\ae{\hbox{\rm a.e.{ }}}

\def\({\Big (}
\def\){\Big )}
\def\[{\Big[}
\def\]{\Big]}

\def\={\buildrel \triangle \over =}

\def\resp{{\it resp. }}

at 9pt 

\def\be{\begin{equation}}
\def\bel{\begin{equation}\label}
\def\ee{\end{equation}}
\def\bea{\begin{eqnarray}}
\def\eea{\end{eqnarray}}
\def\bt{\begin{theorem}}
\def\et{\end{theorem}}
\def\bc{\begin{corollary}}
\def\ec{\end{corollary}}
\def\bl{\begin{lemma}}
\def\el{\end{lemma}}
\def\bp{\begin{proposition}}
\def\ep{\end{proposition}}
\def\br{\begin{remark}}
\def\er{\end{remark}}
\def\ba{\begin{array}}
\def\ea{\end{array}}
\def\bd{\begin{definition}}
\def\ed{\end{definition}}

\newcommand{\inner}[2]{\left\langle#1,#2\right\rangle}
\newcommand{\weakconvergent}[1]{\stackrel{#1}{\relbar\joinrel
\rightharpoonup}}

\newtheorem{lemma}{Lemma}[section]
\newtheorem{remark}{Remark}[section]
\newtheorem{example}{Example}[section]
\newtheorem{theorem}{Theorem}[section]
\newtheorem{corollary}{Corollary}[section]

\newtheorem{definition}{Definition}[section]
\newtheorem{proposition}{Proposition}[section]

\makeatletter
   
   \@addtoreset{equation}{section}
\makeatother

\allowdisplaybreaks

\begin{document}

\title{\bf Second Order Optimality Conditions
for Optimal Control Problems of
Stochastic Evolution Equations}

\author{Qi L\"{u}\thanks{School of Mathematics,
Sichuan University, Chengdu 610064, Sichuan
Province, China. The research of this author is
partially supported by NSF of China under grant
11471231 and the NSFC-CNRS Joint Research
Project under grant 11711530142. {\small\it
E-mail:} {\small\tt lu@scu.edu.cn}.},~~~ Haisen
Zhang\thanks{School of Mathematical Sciences,
Sichuan Normal University, Chengdu 610068,
China.  The research of this author is partially
supported by NSF of China under grants 11401404,
11471231 and 11701470, and the NSF of CQ CSTC
under grants 2015jcyjA00017 and the Advance and
Basic Research Project of Chongqing under grant
cstc2016jcyjA0239. {\small\it E-mail:}
{\small\tt haisenzhang@yeah.net}.}~~~ and~~~ Xu
Zhang\thanks{School of Mathematics, Sichuan
University, Chengdu 610064, Sichuan Province,
China. The research of this author is partially
supported by the NSFC-CNRS Joint Research
Project under grant 11711530142, the PCSIRT
under grant IRT$\_$16R53 and the Chang Jiang
Scholars Program from the Chinese Education
Ministry. {\small\it E-mail:} {\small\tt
zhang$\_$xu@scu.edu.cn}.}}

\date{}

\maketitle

\begin{abstract}
In this paper, we establish some
second order necessary/sufficient optimality
conditions for optimal control problems of 
stochastic evolution equations in infinite dimensions. The control acts
on both the drift and diffusion terms and the
control region is convex. The concepts of
relaxed and $V$-transposition solutions (introduced in our previous works) to operator-valued
backward stochastic evolution equations are employed
to derive these optimality conditions. The
correction part of the second order adjoint
equation, which does not appear in the  (first order) Pontryagin-type
stochastic maximum principle, plays a fundamental role in
our second order optimality conditions.
\end{abstract}

\bs

\no{\bf 2010 Mathematics Subject
Classification}. Primary 93E20; Secondary, 60H07, 60H15.

\bs

\no{\bf Key Words}. Stochastic optimal
control, relaxed transposition solution, $V$-transposition solution, operator-valued backward stochastic evolution equation, second order optimality
condition.

%\newpage

\section{Introduction}\label{s1}

Let $(\O,\cF,\dbF,\dbP)$ be a complete filtered
probability space with the filtration
$\dbF=\{\cF_t\}_{t\ge0}$, on which a
one-dimensional standard Brownian motion
$\{W(t)\}_{t\ge0}$ is defined. Let $T>0$, and
let $X$ be a Banach space with norm
$|\cdot|_{X}$. For any $t\in[0,T]$ and $r\in
[1,\infty)$, denote by $L_{\cF_t}^r(\O;X)$ the
Banach space of all $\cF_t$-measurable random
variables $\xi:\O\to X$ such that
$\mathbb{E}|\xi|_X^r < \infty$, with the
canonical  norm.  Also, denote by
$D_{\dbF}([0,T];L^{r}(\O;X))$ the vector space
of all $X$-valued $\dbF$-adapted processes
$\phi(\cdot)$ such that $\phi(\cdot):[0,T] \to
L^{r}_{\cF_T}(\O;X)$ is c\`adl\`ag, i.e., right
continuous with left limits. Clearly,
$D_{\dbF}([0,T];L^{r}(\O;X))$ is a Banach space
with the following norm
$$
\|\phi(\cd)\|_{D_{\dbF}([0,T];L^{r}(\O;X))} =
\sup_{t\in
[0,T)}\left[\mE|\phi(t)|_X^r\right]^{1/r}.
$$
Denote by $C_{\dbF}([0,T];L^{r}(\O;X))$ the
Banach space of all $X$-valued $\dbF$-adapted
processes $\phi(\cdot)$ such that
$\phi(\cdot):[0,T] \to L^{r}_{\cF_T}(\O;X)$ is
continuous,  with norm inherited from
$D_{\dbF}([0,T];$ $L^{r}(\O;X))$. Fix any
$r_1,r_2,r_3,r_4\in[1,\infty)$. Put
$$
 \ba{ll}
\ds
L^{r_1}_\dbF(\O;L^{r_2}(0,T;X))=\Big\{\f:(0,T)\times\O\to
X\bigm|\f(\cd)\hb{
is $\dbF$-adapted and }\\
\ns\ds\qq\qq\qq\qq\qq\qq\qq

\|\f\|_{L^{r_1}_\dbF(\O;L^{r_2}(0,T;X))}\= \[\dbE\(\int_0^T|\f(t)|_X^{r_2}dt\)^{\frac{r_1}{r_2}}\]^{\frac{1}{r_1}}<\infty\Big\},\\
\ns\ds
L^{r_2}_\dbF(0,T;L^{r_1}(\O;X))=\Big\{\f:(0,T)\times\O\to
X\bigm|\f(\cd)\hb{ is $\dbF$-adapted and
}\\
\ns\ds\qq\qq\qq\qq\qq\qq\qq

\|\f\|_{L^{r_2}_\dbF(0,T;L^{r_1}(\O;X))}\=\[\int_0^T\(\dbE|\f(t)|_X^{r_1}\)^{\frac{r_2}
{r_1}}dt\]^{\frac{1}{r_2}}<\infty\Big\}.
 \ea
$$
Clearly, both $L^{r_1}_\dbF(\O;L^{r_2}(0,T;X))$
and $L^{r_2}_\dbF(0,T;L^{r_1}(\O;X))$ are Banach
spaces with the canonical norms. If $r_1=r_2$,
we simply write the above spaces as
$L^{r_1}_\dbF(0,T;X)$. Let $Y$ be another Banach
space. Denote by $\cL(X; Y)$ the (Banach) space
of all bounded linear operators from $X$ to $Y$,
with the usual operator norm (When $Y=X$, we
simply write $\cL(X)$ instead of $\cL(X; Y)$).
Further, we denote by $
\cL_{pd}\big(L^{r_1}_{\dbF}(0,T;L^{r_2}(\O;X)),\;L^{r_3}_{\dbF}(0,T;L^{r_4}(\O;Y))\big)$
(\resp $
\cL_{pd}\big(X,\;L^{r_3}_{\dbF}(0,T;L^{r_4}(\O;Y))\big)$)
the vector space of all bounded, pointwisely
defined linear operators $\cL$ from
$L^{r_1}_{\dbF}(0,T;L^{r_2}(\O;X))$ (\resp $X$)
to $L^{r_3}_{\dbF}(0,T;L^{r_4}(\O;Y))$, i.e.,
for $\ae (t,\omega)\in (0,T)\times\Omega$, there
exists an $L(t,\o)\in\cL (X;Y)$ verifying that
$\big(\cL \varphi(\cd)\big)(t,\o)=L
(t,\o)\varphi(t,\o), \;\, \forall\;
\varphi(\cd)\in
L^{r_1}_{\dbF}(0,T;L^{r_2}(\O;X))$ (\resp
$\big(\cL x\big)(t,\o)=L (t,\o)x, \;\, \forall\;
x\in X$). Similarly, one can define the spaces
$\cL_{pd}\big(L^{r_2}(\O;X),\;L^{r_3}_{\dbF}(0,T;L^{r_4}(\O;Y))\big)$
and $\cL_{pd}\big(L^{r_2}(\O;X),$
$L^{r_4}(\O;Y)\big)$, etc.

Let $H$  be a separable Hilbert space with the
norm $|\cdot|_{H}$  and the inner product
$\inner{\cdot}{\cdot}_{H}$, and let $A$ be an
unbounded linear operator (with domain $D(A)$ on
$H$), which generates a $C_0$-semigroup
$\{e^{At}\}_{t\geq 0}$. Denote by $A^*$ the
adjoint operator of $A$. Clearly, $D(A)$ is a
Hilbert space with the usual graph norm, and
$A^*$ is the infinitesimal generator of
$\{e^{A^*t}\}_{t\geq 0}$, the adjoint
$C_0$-semigroup of $\{e^{At}\}_{t\geq 0}$. Let
$U$ be a {\it closed convex subset} of another
separable Hilbert space $H_{1}$ (with norm
$|\cdot|_{H_{1}}$  and inner product
$\inner{\cdot}{\cdot}_{H_{1}}$). For any
$\beta\ge 2$, put
$$\cU^{\beta}[0,T] \triangleq \Big\{u(\cd)\in L^{\beta}_\dbF(0,T;H_{1}) \;\Big|\; u(t,\omega)\in U,\; a.e.\; (t,\omega)\in [0,T]\times\Omega\Big\}.$$

Consider the following controlled (forward)
stochastic evolution equation (SEE for short):
\begin{eqnarray}\label{fsystem1}
\left\{
\begin{array}{lll}\ds
dx = \big(Ax +a(t,x,u)\big)dt + b(t,x,u)dW(t) &\mbox{ in }(0,T],\\
\ns\ds x(0)=x_0,
\end{array}
\right.
\end{eqnarray}
where $a$, $b$ are two suitable functions from $[0,T]\times
H\times U$ to $H$, $u\in \cU^{\beta}[0,T]$ and  $x_0\in
L^{\b}_{\cF_0}(\O;H)$. We call $x(\cd)=
x(\cd\,;x_0,u)\in C_\dbF([0,T];L^{\b}(\O;H))$ a
mild solution to \eqref{fsystem1} if
$$
\begin{array}{ll}\ds
x(t)\3n&\ds=e^{At}x_0 + \int_0^t
e^{A(t-s)}a(s,x(s),u(s))ds\\
\ns&\ds\q + \int_0^t
e^{A(t-s)}b(s,x(s),u(s))dW(s),\q
\dbP\mbox{-a.s.,}\q \forall\; t\in[0,T].
\end{array}
$$

Define a cost functional $\cJ(\cdot)$ (for the
control system \eqref{fsystem1}) as follows:
 \bel{jk1}
\cJ(u(\cdot))\triangleq \dbE\Big[\int_0^T
g(t,x(t),u(t))dt + h(x(T))\Big],\qq
u(\cdot)\in \cU^{\beta}[0,T],
 \ee
where $g: [0,T]\times H\times U\to \dbR$ and $h:
H\to \dbR$ are suitably given functions, and $x(\cd)$ is the corresponding
solution to \eqref{fsystem1}.

In this paper we are concerned with the following optimal
control problem for \eqref{fsystem1}:

\ms

\no {\bf Problem (P)} {\it Find a $\bar
u(\cdot)\in \cU^{2}[0,T]$ such that
 \bel{jk2}
\cJ (\bar u(\cdot)) = \inf_{u(\cdot)\in
\cU^{2}[0,T]} \cJ (u(\cdot)).
 \ee
}

Any $\bar u(\cdot)$ satisfying (\ref{jk2}) is
called an {\it optimal control}. The state $\bar
x(\cdot)$ corresponding to $\bar u(\cdot)$  is
called the {\it optimal state}, and
$(\bar x(\cdot),\bar u(\cdot))$ is called an
{\it optimal pair}.

It is one of the most important issues in
optimal control theory  to establish optimality
conditions for optimal controls, which can be employed to distinguish optimal controls from the
other admissible controls. Since the
landmark work in \cite{PC}, first-order
necessary conditions are studied extensively in
the literature for different kinds of control
systems, such as systems governed by ordinary
differential equations (e.g. \cite{PC}), systems
governed by ordinary difference equations (e.g.
\cite{Bol}), systems governed by partial
differential equations (e.g. \cite{LY}), systems
governed by stochastic ordinary differential equations
(e.g. \cite{Peng1, YZ}), systems governed by SEEs (e.g.
\cite{LZ1}), etc.

Similar to the Calculus of Variations (or even the elementary calculus),
in addition to the first-order necessary
conditions, some second order necessary
conditions should be established to distinguish
optimal controls from the candidates which
satisfy the first order necessary conditions,
especially when the optimal controls are
singular, i.e., optimal controls satisfy the
first order necessary conditions trivially. For
instance, when the Hamiltonian corresponding to
optimal controls is equal to a constant in a
subset of the control region or the gradient and
the Hessian (with respect to the control
variable $u$) of the corresponding Hamiltonian
vanish/degenerate. In these cases, the first
order necessary conditions are not enough to
determine the optimal controls. For more
details,  we refer the reader to the
introduction of  \cite{ZZSIAMReview}.

The study of second order necessary conditions
for controlled (deterministic) ordinary differential equations
may date back to the early time of modern
control theory (e.g. \cite{BJ, GK, Goh, Krener}) and attracts
lots of attention until recently (see
\cite{BDP, FT, Knobloch, Lou, Osmolovskii} and the
rich references cited therein). However, as far
as we know, there are merely a few published
papers for second-order necessary conditions for stochastic
optimal control problems  in finite dimensions: 

\begin{itemize}

\item In
\cite{Agayeva, Mahmudov97, Tang}, the main concern
focused on the case that the diffusion term is
independent of the control variable. In \cite{Mahmudov97, Tang}, pointwise second-order
maximum principles for stochastic singular
optimal controls in the sense of Pontryagin-type
maximum principle were established, while in
\cite{Agayeva}, the control system with time
delay was discussed;

\item
When the diffusion terms of
the control systems contain the control
variable, in \cite{BS}, an integral-type
second-order necessary  condition for stochastic
optimal controls was derived under the
assumption that the control region is convex;

\item
Recently, in \cite{ZZ,ZZ2} (see also \cite{ZZSIAMReview}) and \cite{FZhZh},
under some assumptions in terms of the Malliavin
Calculus, the authors established pointwise
second order necessary conditions for stochastic
singular optimal controls with both the convex and the general control
constrains;

\item Very recently, some first and second
order integral type necessary optimality
conditions for stochastic optimal control problems with
state constraints and closed control constraints
were obtained in \cite{FZhZh2,FZhZh3}.
\end{itemize}

The research on the second order sufficient
condition  for optimal controls also has a long
history. It is found that the second order
sufficient condition has important applications
in the sensitivity analysis and the numerical
methods for the optimal control problems. The
corresponding theory for the deterministic cases
has been extensively studied (e.g.
\cite{BDP2, Casas12, Hoehener14, Jacobson70, Lou, Moyer73, Osmolovskii12, Zeidan84}).
However, as far as we know,  \cite{BS} is the
only one reference which contains a sort subsection on the second order sufficient condition for optimal controls of stochastic
control systems in finite dimensions.

To the best of our knowledge, before our work there exists no literature addressed to the second optimality  condition for optimal controls of stochastic
control systems in infinite dimensions.

The main purpose of this paper is to establish
the second order necessary and sufficient
conditions for optimal control problems of SEEs.
In this work, both drift and diffusion terms,
i.e., $a(t,x,u)$ and $b(t,x,u)$, may contain the
control variable $u$, and we assume that the
control region $U$ is convex. The key difference
between \cite{BS, FZhZh, ZZ} and the present work is that
we consider here the SEEs in infinite
dimensions. For such kind of
control systems, the second order adjoint
equation, which is an operator-valued backward
stochastic evolution equation (BSEE for
short), is much more
complex than that in finite dimensions. The
main difficulty to study the well-posedness of
backward stochastic evolution equations is that,
there exists no proper definition
of the It\^{o} integral for operator-valued
stochastic processes (e.g. \cite{van Neerven}).
This leads to some essential
obstacle to obtain the representation of the
correction part of the solutions to such sort of BSEEs. However,  it
can be found in \cite{ZZSIAMReview} that, the
correction part of the second order adjoint
equation plays an indispensable role in the
second order necessary conditions. 

In this
paper, we first employ the notion of relaxed
transposition solution (introduced in \cite{LZ1, LZ}) for the
second order adjoint equations to derive an
integral-type second order necessary condition
for optimal controls. Then, we use the notion of
$V$-transposition solution (introduced in \cite{LZ2}) for the second order
adjoint equations to obtain a pointwise second order
necessary condition. We remark that, quite
different from that in the deterministic setting, there
exist some essential difficulties to derive the
pointwise second-order necessary condition from
an integral-type one when the diffusion term of
the control system depends on the control
variable, \emph{even for the special case of convex
control constraint}. We overcome these
difficulties by some technique developed
in \cite{ZZ}, which is for stochastic control problems in finite dimensions. 

Also, we establish a second
order sufficient condition for optimal controls.
This type of condition essentially ensures that
the cost functional has a quadratic growth
property near an admissible control and hence
ensures the local optimality and uniqueness of
the minimizer. The basic idea comes from the
second order sufficient conditions in
optimization theory.

The rest of this paper is organized as follows:
In Section \ref{s2}, we prove some useful
estimates corresponding to the control system
and present some results for (operator-valued) BSEEs.  Section \ref{s3} is devoted to
establishing the integral type second order
necessary conditions for stochastic optimal
controls. In Section \ref{s4},  we obtain a
pointwise second order necessary optimality condition.  Section \ref{s5}
is addressed to the second order sufficient optimality 
conditions.
Finally,  in
Section \ref{s6} two simple examples are provided to show the applications of the
second order optimality conditions
established in Sections \ref{s4} and \ref{s5}.

Partial results of this paper have been  announced
in \cite{lv2016} without detailed proof.

\ms

\section{Some preliminaries}\label{s2}

Throughout this paper, we assume the following
condition.

\ms

\no{\bf (A1)} {\it Suppose that
$a(\cd,\cd,\cd):[0,T]\times H\times U\to H$ and
$b(\cd,\cd,\cd):[0,T]\times H\times U\to H$ are
two maps satisfying:
\begin{enumerate}[i)]
  \item For any $(x,u)\in
H\times U$, both $a(\cd,x,u):[0,T]\to H$ and
$b(\cd,x,u):[0,T]\to H$ are Lebesgue measurable;
  \item There is a constant $C_L>0$ such that, for a.e. $t\in[0,T]$, any $x,\tilde x\in H$ and any $u,\tilde u\in U$,
\begin{equation}\label{ab0}
\left\{\!\!\!
\begin{array}{ll}\ds
|a(t,x,u) - a(t,\tilde x,\tilde u)|_H+|b(t,x,u) -
b(t,\tilde x,\tilde u)|_H \leq
C_L\big(|x-\tilde x|_H+|u-\tilde u|_{H_{1}}\big),\\
\ns\ds |a(t,0,0)|_H +|b(t,0,0)|_H \leq C_L.
\end{array}
\right.
\end{equation}
\end{enumerate}
}

In the sequel, we shall denote by $C$ a generic
constant, depending on $T$, $A$, $\b$ and $C_L$
(or $F$, $J$ and $K$ to be introduced later),
which may be different from one place to
another. Similar to \cite[Chapter 7]{Prato}, for
any $u(\cdot)\in \cU^{\beta}[0,T]$, it is easy
to show that, under the assumption (A1), the
equation \eqref{fsystem1} is well-posed in the
sense of mild solution and
$$ \|x\|_{C_\dbF([0,T];L^{\b}(\O;H))}\le C\big(1+\|x_{0}\|_{L_{\cF_{0}}^{\b}(\O;H)}+\|u\|_{L^{\beta}_\dbF(\O;L^{2}(0,T;H_{1}))}\big).$$

Also, we need the following condition:

\ms

\no{\bf (A2)} {\it Suppose that
$g(\cd,\cd,\cd):[0,T]\times H\times U\to \dbR$
and $h(\cd):H\to \dbR$ are two functions
satisfying:
\begin{enumerate}[ i)]
  \item For any $(x,u)\in H\times U$,
$g(\cd,x,u):[0,T]\to \dbR$ is Lebesgue
measurable;
  \item  There is a
constant $C_L>0$ such that, for a.e. $t\in[0,T]$, any $x\in H$ and $u\in U$,
\begin{equation}\label{gh}
|g(t,x,u)| +|h(x)|
 \leq C_L(1+|x|_{H}^{2}+|u|_{H_{1}}^{2}).
\end{equation}
\end{enumerate}
}

\ms

Under Assumptions \textbf{(A1)} and
\textbf{(A2)}, the optimal control problem
(\ref{jk2}) (with $\beta\ge2$) is well-defined.

To establish second order necessary conditions,
we need to introduce  further assumptions for
$a(\cd,\cd,\cd)$, $b(\cd,\cd,\cd)$,
$g(\cd,\cd,\cd)$ and $h(\cd)$. To simplify the
notation, for $\varphi= a,b, f$ and $g$, we
denote by $\varphi_{x} (t,x,u)$ and $\varphi_{u}
(t,x,u)$ respectively the first order partial
derivatives of $\varphi$ with respect to $x$ and
$u$ at $(t,x,u)$, by $\varphi_{xx} (t,x,u)$,
$\varphi_{xu} (t,x,u)$ and $\varphi_{uu}
(t,x,u)$ the second order partial derivatives of
$\varphi$ at $(t,x,u)$.

\ms

\no{\bf (A3)} {\it The maps $a(t,\cdot,\cdot)$
and $b(t,\cdot,\cdot)$, and the functional
$g(t,\cdot,\cdot)$ and $h(\cdot)$ are $C^2$ with
respect to $x$ and $u$. Moreover, there exists a
constant $C_L>0$ such that, for a.e. $t\in[0,T]$
and any $(x,u)\in H\times U$,
\begin{equation}\label{ab}
\left\{
\begin{array}{ll}\ds
 \|a_x(t,x,u)\|_{\cL(H)}+\|b_x(t,x,u)\|_{\cL(H)} + \|a_u(t,x,u)\|_{\cL(H_1;H)}+
\|b_u(t,x,u)\|_{\cL(H_1;H)}\leq C_L,\\
\ns\ds
|g_{x}(t,x,u)|_{H} +|g_{u}(t,x,u)|_{H_{1}}+|h_{x}(x)|_{H}
 \leq C_L(1+|x|_{H}+|u|_{H_{1}})\\
\ns\ds \|a_{xx}(t,x,u)\|_{\cL(H\times
H;H)}+\|b_{xx}(t,x,u)\|_{\cL(H\times
H;H)}+\|a_{xu}(t,x,u)\|_{\cL(H\times
H_{1};H)}\\
\ns\ds +\|b_{xu}(t,x,u)\|_{\cL(H\times
H_{1};H)}+\|a_{uu}(t,x,u)\|_{\cL(H_{1}\times
H_{1};H)}+\|b_{uu}(t,x,u)\|_{\cL(H_{1}\times
H_{1};H)}
\leq C_L,\\
\ns\ds
\|g_{xx}(t,x,u)\|_{\cL(H)}+\|g_{xu}(t,x,u)\|_{\cL(H;H_{1})}
+\|g_{uu}(t,x,u)\|_{\cL(H_{1})}+\|h_{xx}(x)\|_{\cL(H)}\leq
C_L.
\end{array}
\right.
\end{equation}}

First, using Assumptions \textbf{(A1)} and
\textbf{(A3)}, we give some estimates for the
control system (\ref{fsystem1}) and its
linearized systems.

Let $\bar u(\cdot)\in \cU^{\beta}[0,T]$ and
$\bar x(\cdot)$ be the corresponding  state of
control system (\ref{fsystem1}). For $\f=a,b$ and $g$,
put
\begin{equation*}
\f_1(t) = \f_x(t,\bar x(t),\bar u(t)),\q \f_2(t)
= \f_u(t,\bar x(t),\bar u(t))
\end{equation*}
and
\begin{equation*}
\f_{11}(t) = \f_{xx}(t,\bar x(t),\bar u(t)),\q
\f_{22}(t) = \f_{uu}(t,\bar x(t),\bar u(t)),\q
\f_{12}(t) = \f_{xu}(t,\bar x(t),\bar u(t)).
\end{equation*}
Let $u(\cdot)\in\cU^{\beta}[0,T]$ be another
admissible control with related state
$x(\cdot)$. Set $\d u(\cdot)=u(\cdot)-\bar
u(\cdot)$ and $\d x(\cdot)= x(\cdot)-\bar
x(\cdot)$. Consider the following first and
second order linearized evolution equations:
\begin{equation}\label{fsystem3.1}
\left\{
\begin{array}{lll}\ds
dy = \big[Ay + a_1(t)y +  a_2(t)\d u
\big]dt + \big[ b_1(t) y + b_2(t)\d u \big]
dW(t) &\mbox{ in
}(0,T],\qq\qq\ \\
\ns\ds y(0)=0;
\end{array}
\right.
\end{equation}
\vspace{+0.5em}
\begin{equation}\label{4.8-eq1}
\left\{
\begin{array}{ll}\ds
dz = \big[Az + a_1(t)z +
a_{11}(t)(y,y) + 2a_{12}(t)(y,\d u) +
a_{22}(t)(\d u,\d u)\big]dt \\
\ns\ds \qq\q + \big[b_1(t)z + b_{11}(t)(y,y)
+2b_{12}(t)(y,\d u) + b_{22}(t)(\d u,\d u)
\big]dW(t) &\mbox{ in
} (0,T],\\
\ns\ds z(0)=0.
\end{array}
\right.
\end{equation}

We have the following estimates.

\begin{lemma}\label{main estmates}
Let {\bf(A1)} and {\bf(A3)} hold. Then, for any $\b\ge2$,
\begin{equation*}
\begin{cases}\ds
\|\delta x\|_{C_\dbF([0,T];L^{\b}(\O;H))}\le
C\|\d u\|_{L^{\b}_\dbF(\Om;L^2(0,T;H_1))},\\
\ns\ds \|y\|_{C_\dbF([0,T];L^{\b}(\O;H))}\le C\|\d u\|_{L^{\b}_\dbF(\Om;L^2(0,T;H_1))},\\
\ns\ds  \|z\|_{C_\dbF([0,T];L^{\b}(\O;H))}\le
C\|\d
u\|^{2}_{L^{2\b}_\dbF(\Om;L^4(0,T;H_1))},\\
\ns\ds \|\delta
x-y\|_{C_\dbF([0,T];L^{\b}(\O;H))} \le C\|\d
u\|^{2}_{L^{2\b}_\dbF(\Om;L^4(0,T;H_1))}. \\
\end{cases}
\end{equation*}
\end{lemma}
{\it Proof}\,: We divide the  proof into two
steps.

{\bf Step 1}. In this step, we prove  the
estimates for $\d x$, $y$ and $z$.

Put
$$
\begin{cases}\ds
\tilde a_1 (t)  \= \int_0^1 a_x(t,\bar x(t) +
\th \d x(t), u(t))d\th,\\
\ns\ds \tilde a_2 (t) \= \int_0^1 a_u(t,\bar
x(t), \bar u(t)+\th\d u(t))d\th,\\
\ns\ds\tilde b_1 (t)  \= \int_0^1 b_x(t,\bar
x(t) + \th \d x(t), u(t))d\th,\\
\ns\ds\tilde b_2 (t) \= \int_0^1 b_u(t,\bar
x(t), \bar u(t)+\th\d u(t))d\th.
\end{cases}
$$

It is easy to see that $\d x(\cdot)$ satisfies
the following SEE:
\begin{equation}\label{fsystem deltax}
\left\{
\begin{array}{lll}\ds
d\d x = \big(A\d x +  \tilde a_1 (t) \d x +
\tilde a_2 (t)\d u  \big)dt + \big( \tilde b_1 (t) \d x +
\tilde b_2 (t)\d u \big)dW(t) &\mbox{ in
}(0,T],\\
\ns\ds \d x(0)=0.
\end{array}
\right.
\end{equation}
Then, by Assumption \textbf{(A3)}, we find that
\begin{equation}\label{estmate deltax}
\begin{array}{ll}\ds
\mE|\d x(t)|^{\b}_H  = \mE\Big| \int_0^t
e^{A(t-s)} \tilde a_1(s) \d x(s) ds + \int_0^t
e^{A(t-s)}
\tilde a_2(s)\d u(s) ds  \\
\ns\ds \hspace{2.4cm} +
\int_0^t e^{A(t-s)}\tilde b_1(s) \d x(s) dW(s) +  \int_0^t e^{A(t-s)}\tilde b_2(s)\d u(s) dW(s)\Big|_H^{\b}\\
\ns\ds \hspace{1.7cm}\leq C \mE\(\Big| \int_0^t
e^{A(t-s)} \tilde a_1(s) \d x(s) ds \Big|_H^{\b}
+  \Big|\int_0^t e^{A(t-s)}\tilde b_1(s) \d x(s) dW(s) \Big|_H^{\b} \\
\ns\ds \hspace{2.4cm}  + \Big| \int_0^t e^{A(t-s)}\tilde a_2(s)\d u(s)ds\Big|_H^{\b} + \Big| \int_0^t e^{A(t-s)}\tilde b_2(s)\d u(s)dW(s)\Big|_H^{\b}\) \\
\ns\ds \hspace{1.7cm} \leq C\[ \int_0^t \mE|\d
x(s)|_H^{\b}ds + \mE \(\int_0^T|\d
u(s)|_{H_1}^{2}ds\)^{\frac{\b}{2}}\].
\end{array}
\end{equation}
It follows from \eqref{estmate deltax} and
Gronwall's inequality that
\begin{equation}\label{estmate1}
\sup_{t\in[0,T]}\mE|\d x(t)|^{\b}_H \leq  C\|\d
u\|^{\b}_{L^{\b}_\dbF(\Om;L^2(0,T;H_1))}.
\end{equation}
In the same way, by \eqref{fsystem3.1} and
Gronwall's inequality we get that
\begin{equation}\label{estmate2}
\sup_{t\in[0,T]}\mE|y(t)|^{\b}_H \leq  C \|\d
u\|^{\b}_{L^{\b}_\dbF(\Om;L^2(0,T;H_1))}.
\end{equation}
Then, by Assumption \textbf{(A3)} and
\eqref{estmate2},
\begin{equation*}
\begin{array}{ll}\ds
\mE|z(t)|^{\b}_H  = \mE\Big| \int_0^t
e^{A(t-s)}\big[
 a_1(s)z(s) +
a_{11}(s)(y(s),y(s)) \\
\ns\ds\hspace{5cm}+ 2a_{12}(s)(y(s),\d u(s)) +
a_{22}(s)(\d u(s),\d u(s))\big]ds  \\
\ns\ds \hspace{2.4cm} +
\int_0^t e^{A(t-s)}\big[b_1(s)z + b_{11}(s)(y(s),y(s))\\
\ns\ds\hspace{5cm}
+2b_{12}(s)(y(s),\d u(s)) + b_{22}(s)(\d u(s),\d u(s))
\big] dW(s)\Big|_H^{\b}\\
\ns\ds \hspace{1.7cm} \leq C\[ \int_0^t
\mE|z(s)|_H^{\b}ds + \mE
\(\int_0^T|y(s)|_{H}^{4} dt\)^{\frac{\b}{2}}
+\mE\(\int_0^T |\d u(s)|_{H_1}^{4}
ds\)^{\frac{\b}{2}}\]\\
\ns\ds \hspace{1.7cm} \leq C\[ \int_0^t
\mE|z(s)|_H^{\b}ds + \mE\int_0^T|y(s)|_{H}^{2\b} dt
+\mE\(\int_0^T |\d u(s)|_{H_1}^{4}
ds\)^{\frac{\b}{2}}\].
\end{array}
\end{equation*}
Therefore,
\begin{equation*}
\sup_{t\in[0,T]}\mE|z(t)|^{\b}_H \leq  C\|\d
u\|^{2\b}_{L^{2\b}_\dbF(\Om;L^4(0,T;H_1))}.
\end{equation*}

\ms

{\bf Step 2}. In this step, we  show that
\begin{equation}\label{estmate3}
\big\|\d x -
y\big\|_{C_\dbF([0,T];L^{\b}(\O;H))}\le C\|\d
u\|^{2}_{L^{2\b}_\dbF(\Om;L^4(0,T;H_1))}.
\end{equation}

Let $r_{1}(\cdot)=\d x(\cdot)-y(\cdot)$. Then
$r_{1}$ solves
\begin{equation}\label{fsystem r1}
\left\{
\begin{array}{lll}\ds
dr_{1} = \big[Ar_{1} +  \tilde a_1 (t) r_{1} +
\big( \tilde a_1 (t)- a_1(t)\big)y  +
\big(  \tilde a_2 (t)-a_2(t)\big)\d u  \big]dt \\
\ns\ds \hspace{1.3cm} + \big[  \tilde b_1 (t) r_{1}
+ \big( \tilde b_1 (t)- b_1(t)\big)y  + \big(
\tilde b_2 (t)-b_2(t)\big)\d u  \big] dW(t) &\mbox{ in
}(0,T],\\
\ns\ds r_{1}(0)=0.
\end{array}
\right.
\end{equation}

By  \textbf{(A3)}, $a_{1}(t, \cdot, \cdot)$,
$a_{2}(t, \cdot, \cdot)$, $b_{1}(t, \cdot,
\cdot)$ and $b_{2}(t, \cdot, \cdot)$ are
Lipschitz on $H\times U$ with respect to $t$
uniformly. Then, it follows from
\eqref{estmate1}, \eqref{estmate2} and
\eqref{fsystem r1} that
$$
\begin{array}{ll}\ds
 \q \mE|r_{1}(t)|^{\b}_{H} \\
\ns\ds = \mE\Big| \int_0^t e^{A(t-s)} \tilde
a_1(s)r_{1}(s)ds
+ \int_0^t e^{A(t-s)} \tilde b_1(s)r_{1}(s) dW(s) \\
\ns\ds\q +\int_0^t e^{A(t-s)}\big( \tilde
a_1(s)-a_1(s)\big) y(s) ds
+\int_0^t e^{A(t-s)}\big( \tilde b_1(s) - b_1(s)\big) y(s) dW(s) \nonumber\\
\ns\ds\q+\int_0^t e^{A(t-s)}\big( \tilde
a_2(s)-a_2(s)\big)\d u(s) ds
+\int_0^t e^{A(t-s)}\big(\tilde b_2(s) - b_2(s)\big)\d u(s)dW(s)\Big|_H^{\b}\\
\ns\ds \leq C\Big\{\mE\int_0^t |r_{1}(s)|_H^{\b}
ds
+\mE\[\int_0^T \big( \|\tilde a_1(s)-a_1(s)\|_{\cL(H)}^{2} + \|\tilde b_1(s)-b_1(s)\|_{\cL(H)}^{2} \big)\cd |y(s)|_H^{2} ds\]^{\frac{\b}{2}} \\
\ns\ds \q+ \mE\[\int_0^T \big(\|\tilde a_2(s) -
a_2(s)\|_{\cL(H_1,H)}^{2}+ \|\tilde b_2(s) -
b_2(s)\|_{\cL(H_1,H)}^{2}\big)\cd|\d u(s)|_{H_1}^{2} ds\]^{\frac{\b}{2}}\Big\}\\
\ns\ds \leq C\[\mE\int_0^t |r_{1}(s)|_H^{\b}
ds+\mE\(\int_0^T \big( |\d x|_{H}^{2} + |\d
u|_{H_1}^{2} \big)\cd |y(s)|_{H}^{2}
ds\)^{\frac{\b}{2}}
+\mE\(\int_0^T|\d u(s)|_{H_1}^{4} ds\)^{\frac{\b}{2}}\]\\
\ns\ds\leq  C\[\mE\int_0^t |r_{1}(s)|_H^{\b}
ds+\mE\(\int_0^T|\d u(s)|_{H_1}^{4}
ds\)^{\frac{\b}{2}}\],
\end{array}
$$
which, together with Gronwall's inequality,
implies that
$$\sup_{t\in[0,T]}\mE|r_{1}(t)|^{\b}_H \leq  C\|\d
u\|^{2\b}_{L^{2\b}_\dbF(\Om;L^4(0,T;H_1))}.$$
This completes the proof of Lemma \ref{main estmates}.
\endpf

\ms

Next, we give a well-posedness result for the
$H$-valued BSEE:
\begin{equation}\label{bsystem1}
\left\{
\begin{array}{lll}
\ds dp(t) = -  A^* p(t) dt + f(t,p(t),q(t))dt + q(t) dW(t) &\mbox{ in }[0,T),\\
\ns\ds p(T) = p_T.
\end{array}
\right.
\end{equation}
Here $p_T \in L_{\cF_T}^{2}(\O;H)$,
$f:[0,T]\times H\times H\times \O \to H$
satisfies
\begin{equation}\label{Lm1}
\left\{
\begin{array}{ll}\ds
f(\cd,0,0)\in
L^{1}_{\dbF}(0,T;L^{2}(\O;H)),\\\ns\ds
|f(t,k_1,k_2)-f(t,\tilde k_1,\tilde k_2)|_H\leq
C_L\big(|k_1-\tilde k_1|_H+|k_2-\tilde k_2|_H \big),\\
\ns\ds\hspace{3.05cm} \ae (t,\o)\in
[0,T]\times\O,\;\; \forall\;k_1,k_2,\tilde
k_1,\tilde k_2\in H.
\end{array}
\right.
\end{equation}
Since neither the usual natural filtration
condition nor the quasi-left continuity is
assumed for the filtration $\dbF$, we cannot
apply the existing results on infinite
dimensional BSEEs (e.g. \cite{Al-H1, HP2, MY, MM})
to obtain the well-posedness of the equation
\eqref{bsystem1}. In what follows we introduce
the concept of the transposition solution to
\eqref{bsystem1} and give the well-posedness
result. To this end, we consider the following
(forward) SEE:
\begin{equation}\label{fsystem2}
\left\{
\begin{array}{lll}\ds
d\f = (A\f+ v_1)ds +  v_2 dW(s) &\mbox{ in }(t,T],\\
\ns\ds \f(t)=\eta,
\end{array}
\right.
\end{equation}
where $t\in[0,T]$, $v_1\in
L^1_{\dbF}(t,T;L^{2}(\O;H))$, $v_2\in
L^2_{\dbF}(t,T;H)$ and $\eta\in
L^{2}_{\cF_t}(\O;H)$ (see \cite[Chapter
6]{Prato} for the well-posedness of
\eqref{fsystem2} in the sense of the mild
solution).

\begin{definition}\label{definition1}
We call $(p(\cdot), q(\cdot)) \in
D_{\dbF}([0,T];L^{2}(\O;H)) \times
L^2_{\dbF}(0,T;H)$  a transposition solution to
(\ref{bsystem1}) if for any $t\in [0,T]$,
$v_1(\cdot)\in L^1_{\dbF}(t,T;L^{2}(\O;H))$,
$v_2(\cdot)\in L^2_{\dbF}(t,T; H)$, $\eta\in
L^{2}_{\cF_t}(\O;H)$ and the corresponding
solution $\f\in C_{\dbF}([t,T];L^{2}(\O;H))$ to
the equation (\ref{fsystem2}), it holds that
\begin{equation*}
\begin{array}{ll}\ds
\q\dbE \big\langle \f(T),p_T\big\rangle_{H}
- \dbE\int_t^T \big\langle \f(s),f(s,p(s),q(s) )\big\rangle_Hds\\
\ns\ds = \dbE \big\langle\eta,p(t)\big\rangle_H
+ \dbE\int_t^T \big\langle
v_1(s),p(s)\big\rangle_H ds + \dbE\int_t^T
\big\langle v_2(s),q(s)\big\rangle_H ds.
\end{array}
\end{equation*}
\end{definition}

\begin{theorem}\cite[Theorem 2.2]{LZ}\label{the1}
Let $p_T \in L^{2}_{\cF_T}(\O;H)$ and
$f(\cd,\cd,\cd)$ satisfy (\ref{Lm1}). Then the
equation \eqref{bsystem1} admits a unique
transposition solution $(p(\cdot), q(\cdot)) \in
D_{\dbF}([0,T];L^{2}(\O; $ $H)) \times
L^2_{\dbF}(0,T;H)$. Furthermore,
\begin{equation*}
\|(p(\cdot), q(\cdot))\|_{
D_{\dbF}([0,T];L^{2}(\O;H)) \times
L^2_{\dbF}(0,T;H)}\leq C\big(\|f(\cd,0,0)\|_{
L^1_{\dbF}(0,T;L^{2}(\O;H))} +\|p_T\|_{
L^{2}_{\cF_T}(\O;H)}\big).
\end{equation*}
\end{theorem}

We also need  the following $\cL(H)$-valued
BSEE\footnote{Throughout this paper, for any
operator-valued process (\resp random variable)
$R$, we denote by $R^*$ its pointwisely dual
operator-valued process (\resp random variable),
e.g., if $R\in L^{r_1}_\dbF(0,T; L^{r_2}(\O;
\cL(H;H_1)))$, then $R^*\in L^{r_1}_\dbF(0,T;
L^{r_2}(\O; \cL(H_1;H)))$, and
$\|R\|_{L^{r_1}_\dbF(0,T; L^{r_2}(\O;
\cL(H;H_1)))}=\|R^*\|_{L^{r_1}_\dbF(0,T;
L^{r_2}(\O; \cL(H_1;H)))}$.}:
\begin{equation}\label{op-bsystem3}
\left\{\3n
\begin{array}{ll}
\ds dP  =  - (A^*  + J^* )P dt  -  P(A + J )dt
-K^*PKdt
 - (K^* Q +  Q K)dt
\\
\ns\ds \hspace{1.21cm}  +   Fdt  +  Q dW(t) &\mbox{ in } [0,T),\\
\ns\ds P(T) = P_T,
\end{array}
\right.
\end{equation}
where
$F\in L^1_\dbF(0,T;L^2(\O;\cL(H)))$,
$P_T\in L^2_{\cF_T}(\O;\cL(H))$, and
$J,K\in L^4_\dbF(0,T; L^\infty(\O;
\cL(H)))$.

To define the solution to \eqref{op-bsystem3},
let us introduce  two SEEs:
\begin{equation}\label{op-fsystem2}
\left\{
\begin{array}{ll}
\ds d\p_1 = (A+J)\p_1ds + u_1ds + K\p_1 dW(s) + v_1 dW(s) &\mbox{ in } (t,T],\\
\ns\ds \p_1(t)=\xi_1
\end{array}
\right.
\end{equation}
and
\begin{equation}\label{op-fsystem3}
\left\{
\begin{array}{ll}
\ds d\p_2 = (A+J)\p_2ds + u_2ds + K\p_2 dW(s) + v_2 dW(s) &\mbox{ in } (t,T],\\
\ns\ds \p_2(t)=\xi_2.
\end{array}
\right.
\end{equation}
Here $\xi_1,\xi_2 \in L^4_{\cF_t}(\O;H)$ and
$u_1,u_2,v_1,v_2\in L^2_\dbF(t,T;L^4(\O;H))$.
Write
\begin{equation*}\label{jshi1}
\begin{array}{ll}\ds
 D_{\dbF,w}([0,T];L^{2}(\O;\cL(H)))\\
\ns\ds\= \Big\{P(\cd,\cd)\;\Big|\;
P(\cd,\cd)\in
\cL_{pd}\big(L^{2}_{\dbF}(0,T;L^{4}(\O;H)),\;L^2_{\dbF}(0,T;L^{\frac{4}{3}}(\O;H))\big),\\
\ns\ds\q \mbox{and for every } t\in[0,T]\hb{ and }\xi\in L^4_{\cF_t}(\O;H),\\
\ns\ds\q P(\cd,\cd)\xi\in
D_{\dbF}([t,T];L^{\frac{4}{3}}(\O;H)) \mbox{ and
}
\|P(\cd,\cd)\xi\|_{D_{\dbF}([t,T];L^{\frac{4}{3}}(\O;H))}
\leq C\|\xi\|_{L^4_{\cF_t}(\O;H)} \Big\}
\end{array}
\end{equation*}
and
\begin{equation*}\label{jshi2}
\3n\begin{array}{ll}\ds
\cQ[0,T]\!\=\!\Big\{\big(Q^{(\cd)},\widehat
Q^{(\cd)}\big)\;\Big|\;\mbox{For any } t\!\in\!
[0,T], \mbox{ both }Q^{(t)}\mbox{ and }\widehat
Q^{(t)}\mbox{ are bounded
linear operators}\\
\ns\ds\hspace{1.7cm}\mbox{ from
}L^4_{\cF_t}(\O;H)\times
L^2_\dbF(t,T;L^4(\O;H))\times
L^2_\dbF(t,T;L^4(\O;H)) \mbox{ to }
L^{2}_\dbF(t,T;L^{\frac{4}{3}}(\O;H))\\
\ns\ds \hspace{1.7cm} \mbox{ and
}Q^{(t)}(0,0,\cd)^*=\widehat
Q^{(t)}(0,0,\cd)\Big\}.
\end{array}
\end{equation*}
\begin{definition}\label{op-definition2x}
We call $\big(P(\cd),(Q^{(\cd)},\widehat
Q^{(\cd)})\big)\in D_{\dbF,w}([0,T];
L^{2}(\O;\cL(H)))\times \cQ[0,T]$
a relaxed transposition solution to
\eqref{op-bsystem3} if for any $t\in
[0,T]$, $\xi_1,\xi_2\in L^4_{\cF_t}(\O;H)$,
$u_1(\cd), u_2(\cd)\in
L^2_{\dbF}(t,T;L^4(\O;H))$ and
$v_1(\cd),v_2(\cd)\in L^2_{\dbF}(t,T;
L^4(\O;H))$, it holds that
 \begin{equation}\label{6.18eq1}
\begin{array}{ll}
\ds \q\mE\big\langle P_T \p_1(T),  \p_2(T) \big\rangle_{H} - \mE \int_t^T
\big\langle
F(s) \p_1(s), \p_2(s) \big\rangle_{H}ds\\
\ns\ds =\mE\big\langle P(t) \xi_1,\xi_2
\big\rangle_{H} + \mE \int_t^T \big\langle
P(s)u_1(s), \p_2(s)\big\rangle_{H}ds + \mE
\int_t^T \big\langle P(s)\p_1(s),
u_2(s)\big\rangle_{H}ds \\
\ns\ds \q  + \mE \int_t^T \big\langle P(s)K
(s)\p_1 (s), v_2 (s)\big\rangle_{H}ds +
\mE
\int_t^T \big\langle  P(s)v_1 (s), K (s)\p_2 (s)+ v_2(s)\big\rangle_{H}ds\\
\ns\ds \q + \mE \int_t^T \big\langle v_1(s),
\widehat
Q^{(t)}(\xi_2,u_2,v_2)(s)\big\rangle_{H}ds+ \mE
\int_t^T \big\langle Q^{(t)}(\xi_1,u_1,v_1)(s),
v_2(s) \big\rangle_{H}ds.
\end{array}
\end{equation}
Here, $\p_1(\cd)$ and $\p_2(\cd)$
solve \eqref{op-fsystem2} and
\eqref{op-fsystem3}, respectively.
\end{definition}

\begin{theorem}\label{OP-th2}
Suppose that  $L^{2}_{\cF_T}(\Om)$  is a separable Banach space. Then, the
equation \eqref{op-bsystem3} admits a unique
relaxed transposition solution
$\big(P(\cd),(Q^{(\cd)},\widehat
Q^{(\cd)})\big) \in  D_{\dbF,w}([0,T];
L^{2}(\O;\cL(H)))$ $\times \cQ[0,T]$. Furthermore,
 $$
\begin{array}{ll}\ds
\q \|P\|_{D_{\dbF,w}([0,T];L^{2}(\Om; \cL(H)))} +
\big\|\big(Q^{(\cd)},\widehat
Q^{(\cd)}\big)\big\|_{\cQ [0,T]}\\[+0.5em]
\ns\ds\leq C\big(
\|F\|_{L^1_\dbF(0,T;\;L^{2}(\Om;\cL(H)))} +
\|P_T\|_{L^{2}_{\cF_T}(\Om;\;\cL(H))}\big).
\end{array}
 $$
\end{theorem}
{\it Proof}\,: See \cite[Chapter 6]{LZ1} or \cite[Section 3]{LZ}.
\endpf

\vspace{0.2cm}

Finally, we introduce the concept of the
$V$-transposition solution to the equation
\eqref{op-bsystem3}. Let  $V$ be a Hilbert space
such that $H\subset V$ and the embedding
operator from $H$ to $V$ is Hilbert-Schmidt.
Denote by $V'$ the dual space of $V$ with
respect to the pivot space $H$. Then we know
that the embedding operator from $V'$ to $H$ is
also Hilbert-Schmidt. Let $X$ and $Y$ be two
Hilbert spaces. Denote by $\cL_2(X;Y)$
($\cL_2(X)$ for $X=Y$) the Hilbert space of all
Hilbert-Schmidt operators from $X$ to $Y$.
\begin{definition}\label{op-definition2}
We call
$$
\begin{array}{ll}\ds
(P(\cd),Q(\cd)) \3n&\ds\in
D_{w,\dbF}([0,T];L^{2}(\Omega;\cL(H)))\times
L^2_{\dbF}(0,T;\cL_2(H;V))
\end{array}
$$
a $V$-transposition solution to
\eqref{op-bsystem3} if for any $t\in [0,T]$,
$\xi_1,\xi_2\in\! L^{4}_{\cF_t}(\Omega;H)$,
$u_1(\cd), u_2(\cd)\in\! L^2_{\dbF}(t,T;$
$L^{4}(\Omega;H))$ and $v_1(\cd), v_2(\cd)\in
L^2_{\dbF}(t,T; L^{4}(\Omega;V'))$, it holds
that
\begin{equation}\label{eq def sol1.1}
\begin{array}{ll}\ds
\q\dbE \big\langle P_T
\p_1(T),\p_2(T)\big\rangle_{H}
 - \dbE\int_t^T \big\langle F(s)\p_1(s),\p_2(s)\big\rangle_{H}ds\\
\ns\ds = \dbE \big\langle
P(t)\xi_1,\xi_2\big\rangle_{H}\! +\!
\dbE\!\int_t^T\!\!\! \big\langle
P(s)u_1(s),\p_2(s)\big\rangle_{H} ds +
\dbE\!\int_t^T\! \big\langle P(s)\p_1(s),u_2(s)\big\rangle_{H} ds \\
\ns\ds \q \!\!+ \dbE\!\int_t^T\! \big\langle
P(s) K(s)\p_1(s), v_2(s)\big\rangle_{H} ds\! +\!
\dbE\!\int_t^T\! \big\langle P(s)v_1\!(s),
K(s)\p_2(s)\!+\!v_2(s)\big\rangle_{H} ds\\
\ns\ds \q \!\! + \dbE\int_t^T \big\langle
v_1(s),Q^*(s)\p_2(s)\big\rangle_{V',V} ds +
\dbE\int_t^T \big\langle
Q(s)\p_1(s),v_2(s)\big\rangle_{V,V'} ds.
\end{array}
\end{equation}
Here, $\p_1(\cd)$ and $\p_2(\cd)$ solve
\eqref{op-fsystem2} and \eqref{op-fsystem3},
respectively.
\end{definition}
Set
$$
\cL_{HV'}\=\Big\{B\in\cL(H)| \mbox{ The restriction
of $B$ on $V'$ belongs to }\cL(V')\Big\}
$$
with the norm
$$
|B|_{\cL_{HV'}} = |B|_{\cL(H)} + |B|_{\cL(V')}.
$$

Let us introduce the following condition:

\vspace{0.2cm}

\no{\bf (A4)} {\it $A$ generates a
$C_0$-semigroup on $V^{'}$ and $J,K\in
L^\infty_\dbF(0,T; \cL_{HV^{'}})$. }

\begin{lemma}\label{op well th0}
{\rm (\cite[Theorem 3.3]{LZ2})} 
Suppose that {\bf (A4)} hold. Then the equation
\eqref{op-bsystem3} admits a unique
$V$-transposition solution $\big(P, Q\big)$.
Furthermore,
\begin{equation}\label{op well th0-eq1}
\begin{array}{ll}\ds
\q\|(P, Q)\|_{
D_{w,\dbF}([0,T];L^{2}(\Omega;\cL(H)))\times
L^2_{\dbF}(0,T;\cL_2(H;V))}\\[+0.5em]
\ns\ds \leq
C\big(\|F\|_{L^1_\dbF(0,T;L^{2}(\Om;\cL(H)))} +
\|P_T\|_{L^{2}_{\cF_T}(\Om;\cL(H))}\big).
\end{array}
\end{equation}
\end{lemma}
%
%%%%%%%%%%%%%%%%%%%%%%%%%%%%%%%%%%%%%%%%%%%%%%%%%%%%%%%%%%%%%%%%%

\section{Integral-type second order necessary conditions} \label{s3}

%%%%%%%%%%%%%%%%%%%%%%%%%%%%%%%%%%%%%%%%%%%%%%%%%%%%%%%%%%%%%%%%

In this section, we  give some integral-type
second order necessary conditions for optimal
controls.

Define
\begin{equation}\label{H}
\begin{array}{ll}\ds
\dbH(t,x,u,k_1,k_2) \= \big\langle k_1,a(t,x,u)  \big\rangle_H + \big\langle k_2, b(t,x,u)  \big\rangle_H - g(t,x,u),\\
\ns\ds \hspace{4cm} (t,x,u,k_1,k_2)\in
[0,T]\times H \times U\times H\times H.
\end{array}
\end{equation}
Let $(\bar x(\cdot), \bar u(\cdot))$ be an
optimal pair, $(p(\cdot),q(\cdot))$ be the
transposition solution of the equation
\eqref{bsystem1}, where $p_T$ and
$f(\cd,\cd,\cd)$ are given by
\begin{equation}\label{zv1}
\left\{
\begin{array}{ll} \ds
p_T =
-h_x\big(\bar x(T)\big),\\\ns
 \ds f(t,k_1,k_2)=-a_x(t,\bar x(t),\bar
u(t))^*k_1 - b_x\big(t,\bar x(t),\bar
u(t)\big)^*k_2 + g_x\big(t,\bar x(t),\bar
u(t)\big).
\end{array}
\right.
\end{equation}
Put
$$
\begin{cases}\ds
\dbH_{x}(t)=\dbH_{x}(t,\bar x(t),\bar
u(t),p(t),q(t)),\\
\ds\dbH_{u}(t)=\dbH_{u}(t,\bar x(t),\bar
u(t),p(t),q(t)),\\
\ds \dbH_{xx}(t) = \dbH_{xx}(t,\bar x(t),\bar
u(t),p(t),q(t)),\\
\ds \dbH_{xu}(t) = \dbH_{xu}(t,\bar x(t),\bar
u(t),p(t),q(t)),\\
\ds \dbH_{uu}(t)= \dbH_{uu}(t,\bar x(t),\bar
u(t),p(t),q(t)).
\end{cases}
$$
Let $(P(\cd),(Q^{(\cd)},\widehat Q^{(\cd)}))$ be
the relaxed transposition solution to the
equation \eqref{op-bsystem3} in which $P_T$,
$J(\cd)$, $K(\cd)$ and $F(\cd)$ are given by
\begin{equation}\label{zv2}
\left\{
\begin{array}{ll} \ds
P_T =
-h_{xx}\big(\bar x(T)\big),\qq\,\;J(t) = a_x(t,\bar x(t),\bar
u(t)),\\
\ns \ds K(t) =b_x(t,\bar x(t),\bar u(t)),\;\;F(t)= -\dbH_{xx}(t).
\end{array}
\right.
\end{equation}

Our main result in this section is as follows.

\begin{theorem}\label{th max}
Assume that $x_0\in L^2_{\cF_0}(\O;H)$ and
$L^2_{\cF_T}(\O)$ is separable. Let
{\bf(A1)}--{\bf(A3)} hold, and let $\bar
u(\cd)\in \cU^{4}[0,T]$ be an optimal control
and $\bar x(\cd)$ be the corresponding optimal
state. Then, for any $u(\cd) \in \cU^{4}[0,T]$
such that
\begin{equation}\label{singular direction}
\mE\int_{0}^{T}\big\langle\dbH_{u}(t), u(t)-\bar
u(t)  \big\rangle_{H_{1}} dt=0,
\end{equation}
the following second order necessary condition holds:
\begin{equation}\label{maxth ine1}
\begin{array}{ll}\ds
\mE\!\int_0^T\! \[ \big\langle
\dbH_{uu}(t)\big(u(t)\!-\!\bar u(t)\big),
u(t)-\bar u(t) \big\rangle_{H_{1}}\! +\!
\big\langle
b_2(t)^{*}P(t)b_2(t)\big(u(t)\!-\!\bar
u(t)\big), u(t)\!-\!\bar u(t)\big\rangle_{H_{1}}
\]dt\\
\ns\ds +2\mE\!\int_0^T\!\!\!
\big\langle\big(\dbH_{xu}(t)\! +\! a_2(t)^*
P(t)\! +\! b_2(t)^*P(t) b_1(t)\big)y(t), u(t)-\bar u(t)\big\rangle_{H_{1}} dt  \\
\ns\ds + \mE\!\int_0^T\! \big\langle\big(
\widehat
Q^{(0)}\!+\!Q^{(0)}\big)\big(0,a_2(t)\big(u(t)-\bar
u(t)\big), b_2(t)\big(u(t)\!-\!\bar
u(t)\big)\big),b_2(t)\big(u(t)\!-\!\bar
u(t)\big)\big) \big\rangle_{H}dt \!\leq\! 0,
\end{array}
\end{equation}
where $y(\cdot)$ is the solution  to the
equation (\ref{fsystem3.1}) corresponding to
$\delta u(\cdot)= u(\cd)-\bar u(\cd)$ and
$(P(\cd),(Q^{(\cd)},\widehat Q^{(\cd)}))$ is
the relaxed transposition solution to the
equation \eqref{op-bsystem3} with the
coefficients given by (\ref{zv2}).
\end{theorem}

{\it Proof}\,: Let us divide the proof into
four steps.

{\bf Step 1}. In this step, we introduce some
notations.

Obviously,  $ \d
u(\cd)=u(\cd)-\bar u(\cd)\in L^4_\dbF(0,T;H_1)$.
Since $U$ is convex, we see that
$$
u^\e(\cdot) = \bar u(\cdot) + \e \d u(\cdot)= (1-\e)\bar u(\cdot) + \e
u(\cdot) \in \cU^{4}[0,T]\subset\cU^{2}[0,T] , \q\forall\;\e \in [0,1].
$$
Denote by $x^\e(\cdot)$ the state process of
\eqref{fsystem1} corresponding to the control
$u^\e(\cdot)$. Let $\d
x^{\e}(\cd)=x^{\e}(\cd)-\bar x(\cd)$ and for
$\psi=a,b,g$, put
$$
\begin{cases}\ds
\tilde{\psi}_{11}^{\e}(t)\=\int_{0}^{1}(1-\theta)\psi_{xx}(t,\bar
x(t)+\theta\delta x^{\e}(t),\bar u(t)+\theta\e
\d u(t))d\theta,\\
\ns\ds
\tilde{\psi}_{12}^{\e}(t)\=\int_{0}^{1}(1-\theta)\psi_{xu}(t,\bar
x(t)+\theta\delta x^{\e}(t),\bar u(t)+\theta\e
\d u(t))d\theta,\\
\ns\ds\tilde{\psi}_{22}^{\e}(t)\=\int_{0}^{1}(1-\theta)\psi_{uu}(t,\bar
x(t)+\theta\delta x^{\e}(t),\bar u(t)+\theta\e
\d u(t))d\theta.\\
\end{cases}
$$
Also, we define
$$
\tilde{h}_{xx}^{\e}(T)\=\int_{0}^{1}(1-\theta)h_{xx}(\bar
x(T)+\theta\delta x^{\e}(T))d\theta.
$$

\vspace{0.2cm}

{\bf Step 2}. It follows from Lemma \ref{main
estmates} that for any $\b\ge 2$,

\begin{equation}\label{max estmate de x and r1}
\begin{cases}\ds
\|\delta x^{\e}\|_{C_\dbF([0,T];L^{\b}(\O;H))}\le C\e\|\d u\|_{L^{\b}_\dbF(\Om;L^2(0,T;H_1))},\;
\\
\ns\ds \|\delta x^{\e}-\e
y\|_{C_\dbF([0,T];L^{\b}(\O;H))} \le C\e^2\|\d
u\|^{2}_{L^{2\b}_\dbF(\Om;L^4(0,T;H_1))}.
\end{cases}
\end{equation}

We claim that there exists a subsequence $\{\e_{n}\}_{n=1}^{\infty}$ such that
\begin{equation}\label{max estmate r2}
\Big\|\delta x^{\e_n}-\e_{n}
y-\frac{\e_{n}^2}{2}z\Big\|_{C_\dbF([0,T];L^{2}(\O;H))}
=o(\e_{n}^2).
\end{equation}
Obviously, $\d x^{\e}$ solves the following SEE:
\begin{equation}\label{exp deltax order two}
\left\{
\begin{array}{ll}
d\delta x^{\e}= \Big[A\delta x^{\e} +a_1(t)\delta x^{\e}
+\e a_{2}(t)\delta u
+\tilde{a}_{11}^{\e}(t)\big(\delta x^{\e},\delta x^{\e}\big)\\[+0.4em]
\qquad\qquad\qq +2\e
\tilde{a}_{12}^{\e}(t)\big(\delta x^{\e},\delta
u\big) +\e^2\tilde{a}_{22}^{\e}(t)\big(\delta
u,\delta u\big)
\Big]dt\\[+0.4em]
\qquad\qquad
+\Big[b_1(t)\delta x^{\e}
+\e b_{2}(t)\delta u
+\tilde{b}_{11}^{\e}(t)\big(\delta x^{\e},\delta x^{\e}\big)\\[+0.4em]
\qquad\qquad\qquad\
+2\e \tilde{b}_{12}^{\e}(t)\big(\delta x^{\e},\delta u\big)
+\e^2 \tilde{b}_{22}^{\e}(t)\big(\delta u,\delta u\big)\Big]dW(t) &\mbox{ in } (0,T],\\[+0.4em]
\delta x^{\e}(0)=0.
\end{array}\right.
\end{equation}
Let $r_{2}^\e(\cd) =\e^{-2}\(\delta
x^{\e}(\cdot)-\e
y(\cdot)-\dfrac{\e^{2}}{2}z(\cdot)\)$. Then
$r_{2}^\e(\cd)$ fulfills
\begin{equation}\label{th3.1-eq3}
\left\{
\begin{array}{ll}\ds
dr_{2}^\e=\Big\{Ar_{2}^\e + a_1(t)r_{2}^\e +\[
a_{11}^\e(t)\(\frac{\d x^\e}{\e},\frac{\d
x^\e}{\e}\)-
\frac{1}{2}a_{11}(t)(y,y)\]\\[+0.5em]
\ns\ds \qq\q
 + \[2\tilde a_{12}^\e(t)\(\frac{\d x^\e}{\e},\d u\) -
a_{12}(t)(y,\d u)\] +
\(\tilde a_{22}^\e(t)-\frac{1}{2}a_{22}(t)\)(\d u,\d u)\Big\}dt \\[+0.5em]
\ns\ds \qq\q\! + \Big\{b_1(t)r_{2}^\e + \[\tilde
b^\e_{11}(t)\(\frac{\d x^\e}{\e},\frac{\d
x^\e}{\e}\) -\frac{1}{2} b_{11}(t)(y,y)\]
\\[+0.5em]
\ns\ds \qq\q\! + \[2\tilde
b_{12}^\e(t)\(\frac{\d x^\e}{\e},\d u\)\! -
b_{12}(t)(y,\d u)\] \!+\! \(\tilde
b_{22}^\e(t)\!-\!\frac{1}{2}b_{22}(t)\)(\d u,\d
u) \Big\}dW(t) \ \ \mbox{ in
}(0,T],\\
\ns\ds r_{2}^\e(0)=0.
\end{array}
\right.
\end{equation}
Put
$$
\begin{array}{ll}\ds
\Psi_{1,\e}(t)\ds=\[ \tilde
a_{11}^\e(t)\(\frac{\d x^\e(t)}{\e},\frac{\d
x^\e(t)}{\e}\)-
\frac{1}{2}a_{11}(t)(y(t),y(t))\]\\
\ns\ds \qq\qq  + \[2\tilde
a_{12}^\e(t)\(\frac{\d x^\e(t)}{\e},\d u(t)\)\!
-\! a_{12}(t)(y(t),\d u(t))\]\! +\! \(\tilde
a_{22}^\e(t)\!-\!\frac{1}{2}a_{22}(t)\)(\d
u(t),\d u(t))
\end{array}
$$
and
$$
\begin{array}{ll}\ds
\Psi_{2,\e}(t)\ds=\[\tilde
b^\e_{11}(t)\(\frac{\d x^\e(t)}{\e},\frac{\d
x^\e(t)}{\e}\) -\frac{1}{2} b_{11}(t)(y(t),y(t))\]
\\
\ns\ds \qq\qq +  \[2\tilde
b_{12}^\e(t)\(\frac{\d x^\e(t)}{\e},\d u(t)\)
\!-\! b_{12}(t)(y(t),\d u(t))\]\! +\! \(\tilde
b_{22}^\e(t)\!-\!\frac{1}{2}b_{22}(t)\)(\d
u(t),\d u(t)).
\end{array}
$$
We have that
\begin{equation}\label{th3.1-eq7}
\begin{array}{ll}\ds
\mE\big|r_2^\e(t)\big|_H^2= \mE\Big| \int_0^t
e^{A(t-s)} a_1(s)r_2^\e(s)ds
+ \int_0^t e^{A(t-s)} b_1(s)r_2^\e(s) dW(s) \\
\ns\ds\qq\qq\qq +\int_0^t
e^{A(t-s)}\Psi_{1,\e}(s) ds +\int_0^t
e^{A(t-s)}\Psi_{2,\e}(s) dW(s) \Big|_H^2
\\
\ns\ds\ \qq\qq \leq C\(\mE\int_0^t |r_2^\e(s)|_H^2
ds + \mE\int_0^t |\Psi_{1,\e}(s)|_H^2 ds +
\mE\int_0^t |\Psi_{2,\e}(s)|_H^2 ds\).
\end{array}
\end{equation}
By (\ref{max estmate de x and r1}), there exists
a  subsequence $\{\e_{n}\}_{n=1}^{\infty}$ such
that $x^{\e_{n}}(\cd)\to \bar x(\cd)$ (in $H$)
a.e. in $\O\times[0,T]$, as $n\to\infty$. Then,
by (\ref{max estmate de x and r1}), Assumption
\textbf{(A3)} and Lebesgue's dominated
convergence theorem, we deduce that
\begin{eqnarray}\label{th3.1-eq8}
&&\lim_{n\to\infty}\mE\int_0^t
|\Psi_{1,\e_{n}}(t)|_H^2 dt\nonumber\\
\!\!\!&\le&\!\!\! \lim_{n\to\infty}\mE\int_0^T\Big| \[
\tilde a_{11}^{\e_{n}}(t)\(\frac{\d
x^{\e_{n}}(t)}{\e_{n}},\frac{\d
x^{\e_{n}}(t)}{\e_{n}}\)-
\frac{1}{2}a_{11}(t)(y(t),y(t))\]\nonumber\\
&& \qq\qq   + \[2\tilde
a_{12}^{\e_{n}}(t)\(\frac{\d
x^{\e_{n}}(t)}{\e_{n}},\d u(t)\) -
a_{12}(t)(y(t),\d u(t))\] \nonumber\\
&& \qq\qq
+\(\tilde a_{22}^{\e_{n}}(t)-\frac{1}{2}a_{22}(t)\)(\d u(t),\d u(t))\Big|_{H}^2dt\nonumber\\
\!\!\!&\le&\!\!\!
C\lim_{n\to\infty}\mE\!\int_0^T\!\!\[\Big|
\tilde a_{11}^{\e_{n}}(t)\(\frac{\d
x{^\e_{n}}(t)}{\e_{n}},\frac{\d
x^\e_{n}(t)}{\e_{n}}\)\!-\! \tilde
a_{11}^{\e_{n}}(t)(y(t),y(t))\Big|_{H}^2\nonumber\\
&&\qq\qq
+\Big\|\tilde
a_{11}^{\e_{n}}(t)\!-\!
\frac{1}{2}a_{11}(t)\Big\|_{\cL(H\times H,\;H)}^2\cdot|y(t)|_{H}^4\nonumber\\
&& \qq\qq +2\Big|\tilde
a_{12}^{\e_{n}}(t)\(\frac{\d
x^{\e_{n}}(t)}{\e_{n}},\d u(t)\) -
\tilde a_{12}^{\e_{n}}(t)(y(t),\d u(t))\Big|_{H}^2\nonumber\\
&& \qq\qq
+\big\|2\tilde a_{12}^{\e_{n}}(t) -
a_{12}(t)\big\|_{\cL(H\times H_{1},\;H)}^2\cdot|y(t)|_{H}^2\cdot|\d u(t)|_{H_{1}}^2\nonumber\\
&& \qq\qq
+
\Big\|\tilde a_{22}^{\e_{n}}(t)-\frac{1}{2}a_{22}(t)\Big\|_{\cL(H_{1}\times H_{1},\;H)}^2\cdot|\d u(t)|_{H_{1}}^4\]dt\nonumber\\
\!\!\!&=&\!\!\!0.
\end{eqnarray}
Similarly,
\begin{equation}\label{th3.1-eq9}
\lim_{n\to\infty}\mE\int_0^t
|\Psi_{2,\e_{n}}(t)|_H^2 dt=0.
\end{equation}
Combining \eqref{th3.1-eq7}, \eqref{th3.1-eq8}
with \eqref{th3.1-eq9}  and using Gronwall's
inequality, we obtain \eqref{max estmate r2}.

\vspace{0.2cm}

{\bf Step 3}. By Taylor's formula, we see that
\begin{equation}\label{th3.1-eq12}
\begin{array}{ll}\ds
\q g(t,x^\e(t),u^\e(t)) - g(t,\bar x(t),
\bar u(t)) \\[+0.4em]
\ns\ds = \big\langle g_1(t),  \delta x^{\e}(t)
\big\rangle_H + \e\big\langle g_2(t),  \d u(t)
\big\rangle_{H_{1}} + \big\langle
\tilde{g}_{11}^{\e}(t)\delta x^{\e}(t), \delta x^{\e}(t) \big\rangle_H \\[+0.4em]
\ns\ds \q +2\e \big\langle
\tilde{g}_{12}^{\e}(t)\delta x^{\e}(t), \d u(t)
\big\rangle_{H_{1}} + \e^2\big\langle
\tilde{g}_{22}^{\e}(t)\d u(t), \d u(t)
\big\rangle_{H_{1}}
\end{array}
\end{equation}
and
\begin{equation}\label{th3.1-eq13}
 h(x^\e(T)) - h(\bar x(T))
= \big\langle h_x(\bar x(T)), \d x^\e(T)
\big\rangle_H + \big\langle
\tilde{h}_{xx}^{\e}(T)\d x^\e(T), \d x^\e(T)
\big\rangle_H.
\end{equation}

Using a similar argument in the proof of
\eqref{max estmate r2}, we can obtain that for
the subsequence $\{\e_{n}\}_{n=1}^{\infty}$ such
that $x^{\e_{n}}(\cd)\to \bar x(\cd)$ (in $H$)
a.e. in $[0,T]\times\O$, as $n\to\infty$,
$$\lim_{n\to \infty}\frac{1}{\e_{n}^2}
\dbE \int_{0}^{T}\Big(\big\langle\tilde{g}_{11}^{\e_{n}}(t)\delta x^{\e_{n}}(t),\delta x^{\e_{n}}(t)\big\rangle_H-\frac{\e_{n}^2}{2}\big\langle g_{11}(t)y(t), y(t)\big\rangle_H\Big)dt=0,$$
$$\lim_{n\to \infty}\frac{1}{\e_{n}^2}
\dbE \int_{0}^{T}\Big(2\big\langle\tilde{g}_{12}^{\e_{n}}(t)\delta x^{\e_{n}}(t),\e_{n} \d u(t)\big\rangle_{H_{1}}-\e_{n}^2 \big\langle g_{12}(t)y(t), \d u(t)\big\rangle_{H_{1}}\Big)dt=0,$$
$$\lim_{n\to \infty}
\dbE
\int_{0}^{T}\Big(\big\langle\tilde{g}_{22}^{\e_{n}}(t)
\d u(t),\d
u(t)\big\rangle_{H_{1}}-\frac{1}{2}\big\langle
g_{22}(t)\d u(t), \d
u(t)\big\rangle_{H_{1}}\Big)dt=0$$ and
$$\lim_{n\to \infty}\frac{1}{\e_{n}^2}
\dbE \Big(\big\langle
\tilde{h}_{xx}^{\e_{n}}(\bar x(T))\delta
x^{\e_{n}}(T),\delta x^{\e_{n}}(T)\big\rangle_H
-\frac{\e_{n}^2}{2}\big\langle h_{xx}(\bar
x(T))y(T), y(T)\big\rangle_H\Big)=0.$$ These,
together with (\ref{max estmate r2}), imply that
\begin{equation}\label{th3.1-eq14}
\begin{array}{ll}\ds
\q \cJ(u^{\e_{n}}) - \cJ(\bar u)\\
\ns\ds = \mE\int_0^T \[\e_{n} \big\langle g_1(t),
y(t)\big\rangle_H + \frac{\e_{n}^2}{2} \big\langle
g_1(t), z(t)\big\rangle_H + \e_{n} \big\langle
g_2(t), \d u(t)\big\rangle_{H_{1}} \\
\ns\ds \q\;\! + \frac{\e_{n}^2}{2}\(\big\langle
g_{11}(t)y(t), y(t)\big\rangle_H + 2\big\langle
g_{12}(t)y(t), \d u(t)\big\rangle_{H_{1}} +
\big\langle g_{22}(t)\d u(t), \d
u(t)\big\rangle_{H_{1}} \)
\]dt \\
\ns\ds \q\;\! + \mE\( \e_{n}\big\langle
h_x(\bar x(T)), y(T) \big\rangle_H \!+\!
\frac{\e_{n}^2}{2}\big\langle h_x(\bar x(T)),
z(T) \big\rangle_H \!+\!
\frac{\e_{n}^2}{2}\big\langle h_{xx}(\bar x(T))
y(T),y(T)\big\rangle_H
\) \!+\! o(\e_{n}^2).
\end{array}
\end{equation}

{\bf Step 4}. By the definition of the
transposition solution to \eqref{bsystem1}, we
have that
\begin{equation}\label{th3.1-eq15}
\begin{array}{ll}\ds
\q\mE\big\langle h_x(\bar x(T)), y(T)
\big\rangle_H\\
\ns\ds =-\mE\int_0^T\(\big\langle
p(t),a_2(t)\d u(t) \big\rangle_H +
\big\langle q(t),b_2(t)\d u(t)
\big\rangle_H +\big\langle g_1(t), y(t)
\big\rangle_H
\)dt
\end{array}
\end{equation}
and
\begin{equation}\label{th3.1-eq16}
\begin{array}{ll}\ds
\q\mE\big\langle h_x(\bar x(T)), z(T)
\big\rangle_H\\
\ns\ds =-\mE\int_0^T\( \big\langle
p(t),a_{11}(t)(y(t),y(t)) \big\rangle_H +
2\big\langle p(t),a_{12}(t)(y(t),\d u(t))
\big\rangle_H \\[+0.4em]
\ns\ds \q + \big\langle p(t),a_{22}(t)(\d
u(t),\d u(t)) \big\rangle_H + \big\langle
q(t),b_{11}(t)(y(t),y(t)) \big\rangle_H \\[+0.4em]
\ns\ds \q  + 2\big\langle q(t),b_{12}(t)(y(t),\d
u(t)) \big\rangle_H + \big\langle
q(t),b_{22}(t)(\d u(t),\d u(t)) \big\rangle_H +
\big\langle g_1(t), z(t) \big\rangle_H \)dt.
\end{array}
\end{equation}
In addition, by the definition of the relaxed transposition
solution to \eqref{op-bsystem3}, we get that
\begin{equation}\label{th3.1-eq17}
\begin{array}{ll}\ds
\q\mE\big\langle h_{xx}(\bar x(T))y(T),y(T)
\big\rangle_H\\[-0.1em]
\ns\ds =-\mE\int_0^T\(\big\langle
P(t)y(t),a_{2}(t)\d u(t) \big\rangle_H +
\big\langle y(t),P(t)a_{2}(t)\d u(t)
\big\rangle_H \\[+0.21em]
\ns\ds \q + \big\langle P(t)b_1(t)y(t),
b_{2}(t)\d u(t) \big\rangle_H + \big\langle
b_1(t)y(t), P(t)b_{2}(t)\d u(t)
\big\rangle_H \\[+0.21em]
\ns\ds \q  + \big\langle P(t)b_2(t)\d u(t),
b_{2}(t)\d u(t) \big\rangle_H + \big\langle
\widehat Q^{(0)}(0,a_2(t)\d
u, b_2(t)\d u)(t), b_2(t)\d u(t) \big\rangle_H \\[+0.21em]
\ns\ds \q + \big\langle Q^{(0)}(0,a_2(t)\d u,
b_2(t)\d u)(t), b_2(t)\d u(t) \big\rangle_H -
\big\langle H_{xx}(t)y(t),y(t)
\big\rangle_H
\)dt.
\end{array}
\end{equation}
Combining \eqref{th3.1-eq14}--\eqref{th3.1-eq17}
with \eqref{singular direction}, we obtain that
\begin{eqnarray*}
0\3n&\le&\3n \frac{\cJ(u^{\e_{n}})-\cJ(\bar u)}{\e_{n}^2}\\[-0.1em]
&=&\3n\! -\mE \int_0^T
\[\frac{1}{\e_{n}}\(\big\langle p(t),a_2(t)\d
u(t) \big\rangle_H + \big\langle q(t),b_2(t)\d
u(t) \big\rangle_H - \big\langle
g_2(t), \d u(t)\big\rangle_{H_{1}} \)\\[+0.1em]
&&\ds\!\!+  \frac{1}{2}\( \big\langle p(t),
a_{22}(t)\big(\d u(t),\d
u(t)\big)\big\rangle_H\!\! +\! \big\langle q(t),
b_{22}(t)\big(\d u(t),\d u(t)\big)\big\rangle_H
\!-\! \big\langle
g_{22}(t)\d u(t), \d u(t)\big\rangle_{H_{1}}\)\\[+0.1em]
&& \ds\!\! + \frac{1}{2}\big\langle P(t)b_2(t)\d
u(t), b_2(t)\d u(t)\big\rangle_H + \(
-\big\langle g_{12}(t)y(t),\d u(t)\big\rangle_{H_{1}}
+ \big\langle p(t),
a_{12}(t)(y,\d u)\big\rangle_H \\[+0.1em]
&&\ds\!\! + \big\langle q(t), b_{12}(t)(y,\d
u)\big\rangle_H + \big\langle a_2(t)^*P(t)y(t),
\d u(t)\big\rangle_{H_{1}} + \big\langle
b_2(t)^*P(t)b_1(t)y(t), \d u(t)\big\rangle_{H_{1}}\\[+0.1em]
&&\ds\!\! + \frac{1}{2}\big\langle \widehat
Q^{(0)}(0,a_2(t)\d u, b_2(t)\d u)(t) +
Q^{(0)}(0,a_2(t)\d u, b_2(t)\d u)(t),
b_2(t)\d u(t) \big\rangle_H \) \]dt + o(1)\\[+0.1em]
\3n&=&\3n\!\! -\mE\!\int_0^T\!\!
\[\frac{1}{\e_{n}} \big\langle \dbH_u(t), \d
u(t) \big\rangle_{H_{1}}\!\!\! +\!
\frac{1}{2}\big\langle \dbH_{uu}(t)\d u(t), \d
u(t) \big\rangle_{H_{1}}\!\!\! +\!
\frac{1}{2}\big\langle b_2(t)^{*}P(t)b_2(t)\d
u(t), \d u(t)\big\rangle_{H_{1}}\!
\]dt\\
\!\!&& \!\!\!\! - \mE\!\int_0^T\!\!\!
\big\langle\big[\dbH_{xu}(t)\! +\! a_2(t)^*
P(t)\! +\! b_2(t)^*P(t) b_1(t)\big]y(t), \d
u(t)\big\rangle_{H_{1}} dt  \\
&& \!\!\!\! - \frac{1}{2}\mE\!\int_0^T\!\!
\big\langle \widehat Q^{(0)}(0,a_2(t)\d u,
b_2(t)\d u)(t) \!+\! Q^{(0)}(0,a_2(t)\d u,
b_2(t)\d u)(t), b_2(t)\d u(t)
\big\rangle_{H}dt \!+\! o(1)\\
\!\!\!&=&\!\!\! -\mE\int_0^T \(
\frac{1}{2}\big\langle \dbH_{uu}(t)\d
u(t), \d u(t) \big\rangle_{H_{1}}\! +\!
\frac{1}{2}\big\langle
b_2(t)^{*}P(t)b_2(t)\d u(t), \d
u(t)\big\rangle_{H_{1}}
\)dt\\
&& \!\!\!\!- \mE\!\int_0^T\!\!\!
\big\langle\big(\dbH_{xu}(t)\! +\! a_2(t)^*
P(t)\! +\! b_2(t)^*P(t) b_1(t)\big)y(t), \d
u(t)\big\rangle_{H_{1}} dt  \\
&& \!\!\!\!-\frac{1}{2}\mE\!\int_0^T\!\!
\big\langle \widehat Q^{(0)}(0,a_2(t)\d u,
b_2(t)\d u)(t)\! +\! Q^{(0)}(0,a_2(t)\d u,
b_2(t)\d u)(t), b_2(t)\d u(t)
\big\rangle_{H}dt\! +\! o(1).
\end{eqnarray*}
Then, letting $n\to \infty$, we finally get
\eqref{maxth ine1}.
\endpf

\vspace{0.3cm}

According to Lemma \ref{op well th0}, to obtain
the well-posedness of \eqref{op-bsystem3} in the
sense of $V$-transposition solution, we only
need the following assumption:

\vspace{0.2cm}

\no{\bf (A5)} {\it $A$ generates a
$C_0$-semigroup on $V'$, $a_x(\cd,\bar
x(\cd),\bar u(\cd)),b_x(\cd,\bar x(\cd),\bar
u(\cd))\!\in\! L^\infty_\dbF(0,\!T; \cL_{HV'})$ and $a_u(\cd,\bar
x(\cd),\bar u(\cd)),b_u(\cd,\bar x(\cd),\bar
u(\cd))\!\in\! L^\infty_\dbF(0,\!T; \cL(H_{1};V'))$. }

\vspace{0.2cm}

Let $(P,Q)$ be the $V$-transposition solution to
BSEE \eqref{op-bsystem3} in which $P_T$,
$J(\cd)$, $K(\cd)$ and $F(\cd)$ are given by
(\ref{zv2}). Put
\begin{equation}\label{S}
\dbS(t)=\dbH_{xu}(t) + a_2(t)^* P(t) +
b_2(t)^*Q(t) +b_2(t)^*P(t)b_1(t).
\end{equation}
The following result holds immediately from
Theorem \ref{th max}.

\begin{corollary}\label{th max v trans}
Let the assumptions in Theorem \ref{th max} and
{\bf (A5)} hold. If $\bar u(\cdot)\in
\cU^{4}[0,T]$, then, for any $u(\cd) \in
\cU^{4}[0,T]$ such that
\begin{equation*}
\mE\int_{0}^{T}\big\langle\dbH_{u}(t), u(t)-\bar u(t)  \big\rangle_{H_{1}} dt=0,
\end{equation*}
the following second order condition holds:
\begin{equation}\label{maxth ine2}
\begin{array}{ll}\ds \mE\int_0^T\[
\big\langle\big(
\dbH_{uu}(t)+b_2(t)^{*}P(t)b_2(t)
\big)\big(u(t)-\bar u(t)\big), u(t)-\bar u(t)
\big\rangle_{H_{1}}\\
\ns\ds \qq\q +2\big\langle \dbS(t)y(t),
u(t)-\bar u(t)\big\rangle_{H_{1}}
\]dt\leq 0.
\end{array}
\end{equation}
\end{corollary}

\section{Pointwise second order necessary conditions} \label{s4}

In this section, we derive the  pointwise second
order necessary condition for  optimal controls
by the integral-type condition (\ref{maxth
ine2}). We assume that $\dbF$ is the natural
filtration generated by $W(\cdot)$. To begin
with, let us introduce some concepts and
technical results which will be used in the rest
of this section.

First, we give the concept of the singular optimal control as follow:

\begin{definition}\label{4.8-def1}
We call $\bar u(\cd)\in \cU^{2}[0,T]$ a
singular optimal control in the classical sense if it is an optimal control and satisfies
\begin{equation}\label{4.8-def1-eq1}
\left\{\!\begin{array}{ll}\ds \!\dbH_u(t,\bar
x(t),\bar u(t),p(t),q(t))=0, \q  a.e.\;
(t,\omega)\in [0,T]\times\Omega,\\[+0.4em]
\ns\ds \lan\big(\dbH_{uu}(t,\bar x(t),\bar
u(t),p(t),q(t)) + b_u(t,\bar x(t),\bar
u(t))^{*}P(t)b_u(t,\bar
x(t),\bar u(t))\big)(v-\bar u(t)), \\
\ns\ds\q v - \bar u(t)\rangle_{H_{1}} = 0,
\qq\forall\ v\in U,\;a.e.\; (t,\omega)\in
[0,T]\times\Omega.
\end{array}
\right.
\end{equation}
\end{definition}

Next, we recall some concepts and results from
Malliavin calculus (see \cite{Nualart} for a
detailed introduction on this topic).

Let $\wt H$ be a separable Hilbert space. We
introduce the Sobolev space $\dbD^{1,2}(\wt H)$
of $\wt H$-valued random variables in the
following way.

Denote by $C_{b}^{\infty}(\dbR^m)$ the set of
$C^{\infty}$-smooth functions with bounded
partial derivatives.  For any $h\in L^2(0,T)$,
write $W(h)=\int_{0}^{T}h(t)dW(t)$. If $F$ is a
smooth $\wt H$-valued random variable of the
form
\begin{equation}\label{1.5-eq1}
F=\sum_{j=1}^n f_j(W(h_{j_{1}}),\cds,W(h_{j_{m}}))\kappa_j
\end{equation}
where $h_{j_{k}}\in L^2(0, T)$, $\kappa_j \in
\wt H$ and $f_j\in
C_{b}^{\infty}(\dbR^{j_{m}})$, $n,j_{m}\in
\dbN$, then the derivative of $F$ is defined as
$$
\cD F=\sum_{j=1}^n\sum_{k=1}^{j_m}
h_{j_{k}}\frac{\pa f_j}{\pa
x_{j_{k}}}(W(h_{j_1}),\cds,W(h_{j_m}))\kappa_j.
$$
Clearly, $\cD F$ is a smooth random variable
with values in $L^2(0,T;\wt H)$. Denote by
$\dbD^{1, 2}(\wt H)$ the completion of the class
of smooth $\wt H$-valued random variables  with
respect to the norm
$$
\|F\|_{\dbD^{1,2}} =\( \mE|F|_{\wt H}^2 +
\mE\int_0^T|\cD_t F|^2_{\wt
H}dt\)^{\frac{1}{2}}.
$$
In particular, given two separable Hilbert
spaces $H_1$ and $H_2$ we can consider $\wt H$ =
$\cL_2(H_1;H_2)$, and in this case, for any $F$
in the space $\dbD^{1,2}(\cL_2(H_1;H_2))$, we
have that $\cD F\in L^2(\Om;L^2(0,T;$
$\cL_2(H_1; H_2))$.

When $\zeta\in \dbD^{1,2}(\wt H)$, the following
Clark-Ocone representation formula holds:
\begin{equation}\label{CO}
\zeta=\mE~\zeta+\int_{0}^{T}\mE(\cD_{s} \zeta\ |\
\cF_{s})dW(s).
\end{equation}
Furthermore, if $\zeta$ is $\cF_{t}$-measurable,
then $\cD_{s}\zeta=0$ for any $s\in(t,T]$.

Write $\dbL^{1,2}(\wt H)$ for the space of
processes $\varphi\in L^{2}([0,T]\times\Omega;
\wt H)$ such that
\begin{enumerate}
  \item[(i)] For a.e. $t\in[0,T]$, $\varphi(t,\cdot)\in \mathbb{D}^{1,2}(\wt H)$;
  \item[(ii)] The function $(s,t,\o)\mapsto\cD_{s}\varphi(t, \omega)$ ($(s,t,\o)\in\ [0,T]\times[0,T]\times\Omega$) admits a measurable version; and
  \item[(iii)] $\displaystyle |||\varphi|||_{1,2}\=\Big(\mE\int_{0}^{T}|\varphi(t)|_{\wt H}^2dt
      +\mE\int_{0}^{T}\int_{0}^{T}|\cD_{s}\varphi(t)|_{\wt H}^2dsdt\Big)^{\frac{1}{2}}<+\infty.$
\end{enumerate}
Denote by $\dbL_{\dbF}^{1,2}(\wt H)$ the set of
all adapted processes in $\dbL^{1,2}(\wt H)$. In
addition, put
\begin{eqnarray*}
& &\dbL_{2^+}^{1,2}(\wt H)\=\Big\{\varphi(\cdot)\in\dbL^{1,2}(\wt H)\Big|\ \exists\ \cD^{+}\varphi(\cdot)\in L^2([0,T]\times\Omega;\wt H)\ \mbox{such that}\\
& &\qquad\qquad\quad
f_{\e}(s)\=\sup_{s<t<(s+\varepsilon)\wedge T}
\mE\big|\cD_{s}\varphi(t)-\cD^{+}\varphi(s)\big|_{\wt H}^2<\infty,\ a.e.\ s\in [0,T],\\
& &\quad\qquad\qquad f_{\e}(\cdot)\ \mbox{is
measurable on }\ [0,T]\ \mbox{for any }
\varepsilon>0,\  \mbox{and}\
\lim_{\varepsilon\to
0^+}\int_{0}^{T}\!\!f_{\e}(s)ds=0\Big\}
\end{eqnarray*}
and
\begin{eqnarray*}
& &\dbL_{2^-}^{1,2}(\wt H)\=\Big\{\varphi(\cdot)\in\dbL^{1,2}(\wt H)\Big|\ \exists\ \cD^{-}\varphi(\cdot)\in L^2([0,T]\times\Omega;\wt H)\ \mbox{such that}\\
& &\qquad\qquad\quad
g_{\e}(s)\=\sup_{(s-\varepsilon)\vee 0<t<s}
\mE\big|\cD_{s}\varphi(t)-\cD^{-}\varphi(s)\big|_{\wt H}^2<\infty,\ a.e.\ s\in [0,T],\\
& &\qquad\qquad\quad g_{\e}(\cdot)\ \mbox{is
measurable on}\ [0,T]\ \mbox{for any }
\varepsilon>0,\ \mbox{and}\ \lim_{\varepsilon\to
0^+}\int_{0}^{T}g_{\e}(s)ds=0\Big\}.
\end{eqnarray*}
Set
$$
\dbL_{2}^{1,2}(\wt H)=\dbL_{2^+}^{1,2}(\wt H)\cap\dbL_{2^-}^{1,2}(\wt H).
$$
For any $\varphi(\cdot)\in \dbL_{2}^{1,2}(\wt
H)$, denote
$\nabla\varphi(\cdot)=\cD^{+}\varphi(\cdot)+\cD^{-}\varphi(\cdot)$.

When $\varphi$ is adapted, $\cD_{s}\varphi(t)=0$
for any $t<s$. In this case,
$\cD^{-}\varphi(\cdot)=0$, and
$\nabla\varphi(\cdot)=\cD^{+}\varphi(\cdot)$.
Denote by $\dbL_{2,\dbF}^{1,2}({\wt H})$ the set
of all adapted processes in $\dbL_{2}^{1,2}({\wt
H})$.

Roughly speaking, an element
$\varphi\in\dbL_{2}^{1,2}(\wt H)$ is a
stochastic process whose  Malliavin derivative
has suitable continuity on some neighbourhood of
$\{(t,t)\ |\ t\in [0,T]\}$. Examples of such
process can be found in \cite{Nualart}.
Especially, if $(s,t)\mapsto \cD_{s}\varphi(t)$
is continuous from $\{(s,t)\big|\ |s-t|<\delta,\
s,t\in [0,T]\}$ (for some $\delta>0$) to
$L_{\cF_{T}}^{2}(\Omega;\wt H)$, then
$\varphi\in\dbL_{2}^{1,2}(\wt H)$ and,
$\cD^{+}\varphi(t)=\cD^{-}\varphi(t)=\cD_{t}\varphi(t)$.

We have the following result.

\begin{lemma}\label{lm4}
Let $\varphi(\cdot)\in \dbL_{2,\dbF}^{1,2}(\wt
H)$. Then, there exists a sequence
$\{\theta_{n}\}_{n=1}^{\infty}$ of positive
numbers such that $\theta_n\to 0^+$ as
$n\to\infty$ and
\begin{equation}\label{lm4-eq1}
\lim_{n\to
\infty}\frac{1}{\theta_{n}^2}\int_{\t}^{\t+\theta_{n}}
\int_{\t}^{t}\mE\big|\cD_{s}\varphi(t)-\nabla
\varphi(s)\big|_{\wt H}^2dsdt =0,\quad a.e.\
\t\in[0,T].
\end{equation}
\end{lemma}
{\em Proof:}\,
For any $\tau, \theta\in [0,\infty)$, we take
the convention that
$$\sup_{t\in [\t,\t+\theta]\cap[0,T]}
\mE\big|D_{\t}\varphi(t)-\nabla\varphi(\t)\big|_{\wt
H}^2=0$$ whenever
$[\t,\t+\theta]\cap[0,T]=\emptyset$. It follows
from the definition of $\dbL_{2,\dbF}^{1,2}(\wt
H)$ that
\begin{eqnarray*}
&&\q\lim_{\theta\to
0^+}\frac{1}{\theta^2}\int_{0}^{T}\int_{\t}^{\t+\theta}
\int_{\t}^{t}\mE\big|\cD_{s}\varphi(t)-\nabla \varphi(s)\big|_{\wt H}^2dsdtd\t\\
&&\ds= \lim_{\theta\to
0^+}\frac{1}{\theta^2}\int_{0}^{T}\int_{\t}^{\t+\theta}
\int_{s}^{\t+\theta}\mE\big|\cD_{s}\varphi(t)-\nabla\varphi(s)\big|_{\wt H}^2dtdsd\t\\
&&\ds\leq\lim_{\theta\to
0^+}\frac{1}{\theta}\int_{0}^{T}\int_{\t}^{\t+\theta}
\Big[\sup_{t\in [s,s+\theta]\cap[0,T]}\mE
\big|\cD_{s}\varphi(t)-\nabla\varphi(s)\big|_{\wt H}^2\Big]dsd\t\\
&&\ds\leq\lim_{\theta\to
0^+}\frac{1}{\theta}\int_{0}^{T}\int_{0}^{\theta}
\Big[\sup_{t\in [s+\t,s+\t+\theta]\cap[0,T]}
\mE\big|\cD_{s+\t}\varphi(t)
-\nabla \varphi(s+\t)\big|_{\wt H}^2\Big]dsd\t\\
&&\ds\leq\lim_{\theta\to
0^+}\frac{1}{\theta}\int_{0}^{\theta}\int_{0}^{T}
\Big[\sup_{t\in [s+\t,s+\t+\theta]\cap[0,T]}
\mE\big|\cD_{s+\t}\varphi(t)-\nabla\varphi(s+\t)\big|_{\wt H}^2\Big]d\t ds \\
&&\ds\leq\lim_{\theta\to
0^+}\frac{1}{\theta}\int_{0}^{\theta}\int_{s}^{T}
\Big[\sup_{t\in [\t,\t+\theta]\cap[0,T]}
\mE\big|\cD_{\t}\varphi(t)-\nabla\varphi(\t)\big|_{\wt H}^2\Big]d\t ds \\
&&\ds\leq\lim_{\theta\to
0^+}\frac{1}{\theta}\int_{0}^{\theta}\int_{0}^{T}
\Big[\sup_{t\in [\t,\t+\theta]\cap[0,T]}
\mE\big|\cD_{\t}\varphi(t)-\nabla\varphi(\t)\big|_{\wt H}^2\Big]d\t ds \\
&&\ds\leq\lim_{\theta\to 0^+}\int_{0}^{T}
\Big[\sup_{t\in [\t,\t+\theta]\cap[0,T]}
\mE\big|\cD_{\t}\varphi(t)-\nabla\varphi(\t)\big|_{\wt H}^2\Big]d\t  \\
&&\ds=  0,
\end{eqnarray*}
which implies \eqref{lm4-eq1}.
\endpf

The following results will be frequently used in
the proof of the main results in this section.

\begin{lemma}\label{lm3}
Let $\phi(\cdot),\ \psi(\cdot)\in
L_{\dbF}^{2}(0,T;H)$. Then, for a.e. $\t\in
[0,T)$,
\begin{equation}\label{lm3-eq1}
\lim_{\theta\to
0^+}\frac{1}{\theta^2}\dbE\int_{\t}^{\t+\theta}\Big\langle\phi(\t),
\int_{\t}^{t} e^{A(t-s)}\psi(s)ds\Big\rangle_{H}
dt
=\frac{1}{2}\dbE\langle
\phi(\t),\psi(\t)\rangle_{H},
\end{equation}
\begin{equation}\label{lm3-eq2}
\lim_{\theta\to
0^+}\frac{1}{\theta^2}\dbE\int_{\t}^{\t+\theta}
\Big\langle\phi(t), \int_{\t}^{t}e^{A(t-s)}
\psi(s)ds \Big\rangle_{H} dt
=\frac{1}{2}\dbE\langle\phi(\t),\psi(\t)\rangle_{H}.
\end{equation}
\end{lemma}

{\it Proof}\,: The equality \eqref{lm3-eq1} is a
corollary of the Lebesgue differentiation
theorem. Now, we prove \eqref{lm3-eq2}. For any
$\t\in [0,T)$, let $\theta>0$ and $\t+\theta<T$.
It follows from the Lebesgue differentiation
theorem that
$$
\lim_{\theta\to
0^+}\frac{1}{\theta}\int_{\t}^{\t+\theta}
\mE\big|\phi(t)-\phi(\t)\big|_{H}^2dt = 0,\ \ \
a.e.\ \t\in [0,T),
$$
and
$$
\lim_{\theta\to
0^+}\frac{1}{\theta^{2}}\mE\int_{\t}^{\t+\theta}
\int_{\t}^{t}\big|e^{A(t-s)}\psi(s)\big|_{H}^2dsdt
= \frac{1}{2}\mE\big|\psi(\t)\big|_{H}^2,\ \ \
a.e.\ \t\in [0,T).
$$
Therefore,
\begin{equation}\label{1.5-eq2}
\begin{array}{ll}\ds
\q\lim_{\theta\to
0^+}\Big|\frac{1}{\theta^2}\mE\int_{\t}^{\t+\theta}
\Big\langle \phi(t)-\phi(\t), \int_{\t}^{t} e^{A(t-s)}\psi(s)ds\Big\rangle_{H} dt\Big|\\
\ns\ds\leq \lim_{\theta\to
0^+}\frac{1}{\theta^2}
\Big[\int_{\t}^{\t+\theta}
\mE\big|\phi(t)-\phi(\t)\big|^2_{H}dt\Big]^{\frac{1}{2}}
\Big[\int_{\t}^{\t+\theta}
(t-\t)\int_{\t}^{t}\mE\big|e^{A(t-s)}\psi(s)
\big|_{H}^2ds
dt\Big]^{\frac{1}{2}}\\
\ns\ds\leq\lim_{\theta\to
0^+}\frac{1}{\theta^{\frac{3}{2}}}\Big[\int_{\t}^{\t+\theta}
\mE\big|\phi(t)-\phi(\t)\big|_{H}^2dt\Big]^{\frac{1}{2}}
\Big[\int_{\t}^{\t+\theta} \int_{\t}^{t}
\mE\big|e^{A(t-s)}\psi(s)\big|_{H}^2ds
dt\Big]^{\frac{1}{2}}
 \\
\ns\ds =  0,\ \qquad a.e.\ \t\in [0,T).
\end{array}
\end{equation}
From \eqref{lm3-eq1}  and \eqref{1.5-eq2}, we
obtain \eqref{lm3-eq2}. This completes the proof
of Lemma \ref{lm3}.
\endpf

\vspace{0.2cm}

Now, we assume that

\medskip

\no{\bf (A6)} \vspace{-0.4cm}  {\it
$$\bar{u}(\cdot)\in\dbL_{2,\dbF}^{1,2}(H_{1}),\;
\dbS(\cdot)^*\in
\dbL_{2,\dbF}^{1,2}(\cL_2(H_{1};H))\cap
L^{\infty}([0,T]\times\Omega;\cL_2(H_{1};H)),$$
and
$$\cD_{\cdot}\dbS(\cdot)^*\in
L^{2}(0,T;L^{\infty}([0,T]\times\Omega;\cL_2(H_{1};H))).$$
}
\begin{remark}
{\bf (A6)} is a restriction on the regularity of
optimal controls. We believe that it is a
technical condition.  However, we do not know
how to get rid of it now.
\end{remark}

\begin{remark}
We can replace {\bf (A6)} by the following
assumption:

\ms

\no{\bf (A6')}   {\it \vspace{-0.65cm}
$$\bar{u}(\cdot)\in\dbL_{2,\dbF}^{1,2}(H_{1}),\
\dbS(\cdot)^{*} \in
\dbL_{2,\dbF}^{1,2}(\cL_2(H_1;V'))\cap
L^{\infty}([0,T]\times\Omega;\cL_2(H_1;V')),
$$
$$\cD_{\cdot}\dbS(\cdot)^*\in
L^{2}(0,T;L^{\infty}([0,T]\times\Omega;\cL_2(H_{1};V'))),$$
and
$$
a_u(\cd,\bar x(\cd),\bar
u(\cd)),b_u(\cd,\bar x(\cd),\bar
u(\cd))\in L^\infty_\dbF(0,T;\cL(H_{1};V')).
$$
}

In {\bf (A6')}, we relax the restriction of the
regularity on $H$ by assuming that $a_u,b_u$ can
map $H_{1}$ into a more regular space $V'$.
\end{remark}

By Assumption {\bf (A6)}, for any $v\in U$,
$\dbS(t)^{*}(v-\bar{u}(t))\in
\dbL_{2,\dbF}^{1,2}(H)$ and
\begin{equation}\label{12.20-eq1}
\dbS(t)^{*}(v-\bar{u}(t)) =\mE
~\Big[\dbS(t)^{*}(v-\bar{u}(t))\Big]
+\int_{0}^{t}
\mE\[\cD_{s}\(\dbS(t)^{*}(v-\bar{u}(t))\)\ \Big|
\cF_{s}\]  dW(s),\quad \dbP\mbox{-}a.s.
\end{equation}

Now we are about to give our main result,  the
pointwise second order necessary condition for
singular optimal controls.
When the optimal control $\bar u$ is singular in
the sense of  Definition \ref{4.8-def1},  the
following result is an immediate consequence of
Corollary \ref{th max v trans}.
\begin{corollary}\label{th max1}
Assume that $x_0\in L^2_{\cF_0}(\O;H)$. Let
Assumptions  {\bf(A1)}--\textbf{\bf(A3)},
{\bf(A5)} hold, and let $\bar
u(\cd)\in\cU^{4}[0,T]$ be a singular optimal
control and $\bar x(\cd)$ be the corresponding
optimal state. Then, for any $u(\cd) \in
\cU^{4}[0,T]$,
\begin{equation}\label{maxth ine1.1}
\begin{array}{ll}\ds
\mE \int_0^T \big\langle y(t), \dbS(t)^*\big(u(t)
-\bar u(t)\big)\big\rangle_{H} dt  \leq 0.
\end{array}
\end{equation}
\end{corollary}

Using \eqref{12.20-eq1} and  \eqref{maxth ine1.1} ,
we have the following pointwise second-order
necessary condition for singular optimal
controls.
\begin{theorem}\label{th max p1}
Let Assumptions {\bf(A1)}--{\bf(A3)} and
{\bf(A5)}--{\bf(A6)} hold. If $\bar{u}(\cdot)\in
\cU^{4}[0,T]$ is a singular optimal control in
the classical sense, then for a.e. $\t\in
[0,T]$, it holds that
\begin{equation}\label{th max p1-eq1}
\begin{array}{ll}\ds
\langle a_{2}(\t)(v-\bar{u}(\t)),\dbS(\t)^{*}(v-\bar{u}(\t))\rangle_{H }
+ \langle b_{2}(\t)(v-\bar{u}(\t)),\nabla
\dbS(\t)^{*}(v-\bar{u}(\t))\rangle_{H }\\[+0.4em]
\ns\ds
 -
\langle b_{2}(\t)(v-\bar{u}(\t)),\dbS(\t)^{*}\nabla\bar{u}(\t)\rangle_{H }
 \le 0, \quad \forall \ v\in U,
\ \dbP\mbox{-}a.s.
\end{array}
\end{equation}
\end{theorem}
{\em Proof:}\, Since $W(\cdot)$ is a continuous
stochastic process, $\cF_{t}$ is countably
generated for any $t\in[0,T]$. Hence, one can
find a sequence $\{F_{l}\}_{l=1}^{\infty}\subset
\cF_{t}$ such that for any $F\in \cF_{t}$, there
exists a subsequence
$\{F_{l_{n}}\}_{n=1}^{\infty}\subset
\{F_{l}\}_{l=1}^{\infty}$ such that
$\ds\lim_{n\to \infty} \dbP\big((F\setminus
F_{l_{n}})\bigcup (F_{l_{n}}\setminus
F)\big)=0$. $\cF_{t}$ is also said to be
generated by the sequence
$\{F_{l}\}_{l=1}^{\infty}$.

Denote by $\{t_{i}\}_{i=1}^{\infty}$ the
sequence constituted by all rational numbers in
$[0,T)$, by $\{v^{k}\}_{k=1}^{\infty}$ a dense
subset of $U$. As in \cite{Haussmann}, we choose
$\{F_{ij}\}_{j=1}^{\infty}(\subset \cF_{t_i})$
to be a sequence generating $\cF_{t_i}$ (for
each $i\in \dbN$). Fix $i,j,k\in \dbN$
arbitrarily. For any $\t\in [t_{i},T)$ and
$\theta\in(0,T-\t)$, write
$E_{\theta}^{i}=[\t,\t+\theta)$, and define
$$
u_{ij}^{k,\theta}(t,\omega)=\left\{
\begin{array}{l}
v^{k}, \qquad\qquad\quad (t,\omega)\in E_{\theta}^{i}\times F_{ij},\\[+0.4em]
\bar{u}(t,\omega), \qquad \quad\, (t,\omega)\in \big([0,T]\times \Omega\big) \setminus  \big(E_{\theta}^{i}\times F_{ij}\big).\\
\end{array}\right.
$$
Clearly, $u_{ij}^{k,\theta}(\cdot)\in \cU^{4}[0,T]$ and
$$
u_{ij}^{k,\theta}(t,\omega)-\bar{u}(t,\omega)
=\big(v^{k}-\bar{u}(t,\omega)\big)\chi_{F_{ij}}(\omega)\chi_{E_{\theta}^{i}}(t),\quad
(t,\omega)\in [0,T]\times \Omega.
$$
Then, substituting $u(\cdot)$ by $u_{ij}^{k,\theta}(\cd)$ in
\eqref{maxth ine1.1}, we obtain that
\begin{equation}\label{12.20-eq35}
\dbE\int_{\t}^{\t+\theta} \big\langle
y_{ij}^{k,\th}(t),\dbS(t)^*\big(v^{k}-\bar{u}(t)\big)\big\rangle_{H}\chi_{F_{ij}}(\omega)dt\le
0,
\end{equation}
where $y_{ij}^{k,\th}(\cdot)$ is the solution to
the equation (\ref{fsystem3.1}) with $u(\cdot)$
replaced by $u_{ij}^{k,\theta}(\cdot)$. Note
that $y_{ij}^{k,\th}(\cdot)$ is the mild
solution to the linear evolution equation
(\ref{fsystem3.1}), i.e.,
\begin{equation}\label{12.20-eq8}
\begin{array}{ll}\ds
y_{ij}^{k,\th}(t)\3n&\ds=
\int_{0}^{t}e^{A(t-s)}
\[a_1(s)y_{ij}^{k,\th}(s) +  a_2(s)\big(v^{k}-\bar{u}(s)\big)\chi_{E_{\theta}^{i}}(s)\chi_{F_{ij}}(\omega)\]ds\\
\ns&\ds \q + \int_{0}^{t}e^{A(t-s)}\[
b_1(s)y_{ij}^{k,\th}(s) +
b_2(s)\big(v^{k}-\bar{u}(s)\big)\chi_{E_{\theta}^{i}}(s)\chi_{F_{ij}}(\omega)
\]dW(s),\\
\ns&\ds \hspace{9cm} \dbP\mbox{-a.s.,} \
\forall\ t\in[0,T].
\end{array}
\end{equation}
Substituting (\ref{12.20-eq8}) into
(\ref{12.20-eq35}) and recalling that
$y_{ij}^{k,\th}(t)=0$ for any $t\in [0,\t)$, we
have
\begin{equation}\label{12.20-eq9}
\begin{array}{ll}\ds
0\3n&\ds\ge
\frac{1}{\theta^2}\dbE\int_{\t}^{\t+\theta}
\Big\langle \int_{\t}^{t}e^{A(t-s)}
\[a_1(s)y_{ij}^{k,\th}(s) \\
\ns&\ds \qquad\qquad\qquad\qquad\qquad\qquad
+a_2(s)\big(v^{k}-\bar{u}(s)\big)\chi_{F_{ij}}(\omega)\]ds,
 \dbS(t)^*\big(v^{k}-\bar{u}(t)\big)\Big\rangle_{H}\chi_{F_{ij}}(\omega)dt \\
\ns&\ds+\frac{1}{\theta^2}
\dbE\int_{\t}^{\t+\theta} \Big\langle
\int_{\t}^{t}
e^{A(t-s)}\[ b_1(s)y_{ij}^{k,\th}(s)\\
\ns&\ds\qquad\qquad\qquad\qquad\qquad
+ b_2(s)\big(v^{k}-\bar{u}(s)\big)\chi_{F_{ij}}(\omega)
\]dW(s),
 \dbS(t)^*\big(v^{k}-\bar{u}(t)\big)\Big\rangle_{H}\chi_{F_{ij}}(\omega)dt.
\end{array}
\end{equation}
By \textbf{(A3)}, using  Gronwall's inequality
and Burkholder-Davis-Gundy's inequality, it is
easy to prove that
\begin{equation}\label{est of x1ijk}
\sup_{t\in[0,T]}\dbE|y_{ij}^{k,\th}(t)|_{H}^{2}\le
C\dbE\int^{T}_{0}|v^{k}-\bar u(s)|_{H_{1}}^{2}\chi_{E_{\theta}^{i}}(s)\chi_{F_{ij}}(\o)ds.
\end{equation}
Consequently, for a.e. $\t\in [t_{i},T)$,
\begin{equation}\label{12.20-eq10add}
\begin{array}{ll}\ds
\q \frac{1}{\theta^2}
\Big|\dbE\int_{\t}^{\t+\theta} \Big\langle
\int_{\t}^{t}e^{A(t-s)}
a_1(s)y_{ij}^{k,\th}(s)ds,
 \dbS(t)^*\big(v^{k}-\bar{u}(t)\big)\Big\rangle_{H}\chi_{F_{ij}}(\omega)dt\Big|\\
\ns\ds\le\frac{1}{\theta^2} \Big(\dbE~
\int_{\t}^{\t+\theta}\Big|\int_{\t}^{t}e^{A(t-s)}
a_1(s)y_{ij}^{k,\th}(s)ds\Big|_{H}^2dt\Big)^{\frac{1}{2}}
\Big(\dbE~ \int_{\t}^{\t+\theta}\Big|\dbS(t)^{*}\big(v^{k}-\bar{u}(t)\big) \Big|_H^2dt\Big)^{\frac{1}{2}}\\
\ns\ds\le\frac{1}{\theta^2} \Big(\dbE~
\int_{\t}^{\t+\theta}(t\!-\!\t)\int_{\t}^{t}\!\big|e^{A(t-s)}
a_1(s)y_{ij}^{k,\th}(s)\big|_{H}^2dsdt\Big)^{\frac{1}{2}}
 \Big(\dbE~ \int_{\t}^{\t+\theta}\Big|\dbS(t)^{*}\big(v^{k}\!-\!\bar{u}(t)\big) \Big|_H^2dt\Big)^{\frac{1}{2}}\\
\ns\ds\le\frac{C}{\theta^{\frac{1}{2}}}
\Big(\sup_{t\in[0,T]}\dbE|y_{ij}^{k,\th}(t)|_{H}^{2}\Big)^{\frac{1}{2}}
\Big(\dbE~
\int_{\t}^{\t+\theta}\Big|\dbS(t)^{*}\big(v^{k}-\bar{u}(t)\big)
\Big|_H^2dt\Big)^{\frac{1}{2}}\to 0, \q
\theta\to 0^+.
\end{array}
\end{equation}
Next, by Lemma \ref{lm3},  for a.e. $\t\in [t_{i},T)$,
\begin{equation}\label{12.20-eq10}
\begin{array}{ll}\ds
\q \lim_{\theta\to
0^+}\frac{1}{\theta^2}\dbE\int_{\t}^{\t+\theta}
\Big\langle \int_{\t}^{t}\Big[e^{A(t-s)}
a_2(s)\\
\ns\ds\qq\qq\qq\qq\qq\qq
\big(v^{k}-\bar{u}(s)\big)\chi_{F_{ij}}(\omega)\Big]ds,
 \dbS(t)^*\big(v^{k}-\bar{u}(t)\big)\Big\rangle_{H}\chi_{F_{ij}}(\omega)dt\\
\ns\ds= \frac{1}{2}\dbE~\big(\langle
a_{2}(\t)(v^{k}-\bar{u}(\t)),
\dbS(\t)^*\big(v^{k}-\bar{u}(\t)\big)\rangle_{H}\chi_{F_{ij}}(\omega)\big).
\end{array}
\end{equation}
Therefore, by (\ref{12.20-eq10add}) and (\ref{12.20-eq10}), we have already proved that
\begin{equation}\label{12.20-limit part1}
\begin{array}{ll}\ds
\lim_{\theta\to
0^+}\frac{1}{\theta^2}\dbE\int_{\t}^{\t+\theta}
\Big\langle \int_{\t}^{t}e^{A(t-s)}
\[a_1(s)y_{ij}^{k,\th}(s) \\
\ns\ds \qquad\qquad\qquad\qquad\qquad
+a_2(s)\big(v^{k}-\bar{u}(s)\big)\chi_{F_{ij}}(\omega)\]ds,
 \dbS(t)^*\big(v^{k}-\bar{u}(t)\big)\Big\rangle_{H}\chi_{F_{ij}}(\omega)dt \\
\ns\ds= \frac{1}{2}\dbE~\big(\langle
a_{2}(\t)(v^{k}-\bar{u}(\t)),
 \dbS(\t)^*\big(v^{k}-\bar{u}(\t)\big)\rangle_{H}\chi_{F_{ij}}(\omega)\big), \qq a.e.\ \t\in[t_{i},T).
\end{array}
\end{equation}

On the other hand, by Assumption \textbf{(A6)}
and (\ref{12.20-eq1}),
\begin{eqnarray}\label{12.20-eq11}
&&\ds \dbE\int_{\t}^{\t+\theta} \Big\langle
\int_{\t}^{t}
e^{A(t-s)}\[ b_1(s)y_{ij}^{k,\th}(s) \nonumber\\
&&\ds\qquad\qquad\qquad\qquad\q +
b_2(s)\big(v^{k}-\bar{u}(s)\big)\chi_{F_{ij}}(\omega)
\]dW(s),
\dbS(t)^*\big(v^{k}-\bar{u}(t)\big)\Big\rangle_{H}\chi_{F_{ij}}(\omega)dt \nonumber\\
&&\ds=\int_{\t}^{\t+\theta}
\dbE~\Big\{\Big\langle\int_{\t}^{t}e^{A(t-s)}\[
b_1(s)y_{ij}^{k,\th}(s)+
b_2(s)\big(v^{k}-\bar{u}(s)\big)\chi_{F_{ij}}(\omega)
\]dW(s),\nonumber\\
&&\ds
\qquad\qquad\qquad\qquad\qquad\qquad\qquad\qquad\qquad
\dbE~\big[
\dbS(t)^*\big(v^{k}-\bar{u}(t)\big)\big]
\Big\rangle_{H}\chi_{F_{ij}}(\omega)\Big\}dt \\
&&\ds \qquad+\int_{\t}^{\t+\theta}
\dbE~\Big\{\Big\langle\int_{\t}^{t}e^{A(t-s)}\[
b_1(s)y_{ij}^{k,\th}(s)+
b_2(s)\big(v^{k}-\bar{u}(s)\big)\chi_{F_{ij}}(\omega)
\]dW(s),\nonumber\\
&&\ds \qquad\qquad\qquad\qquad\quad\ \
\int_{0}^{t}\dbE~\Big[\cD_{s}
\big(\dbS(t)^{*}(v^{k}-\bar{u}(t))\big)\ \Big|\
\cF_{s}\Big]dW(s)
\Big\rangle_{H}\chi_{F_{ij}}(\omega)\Big\}dt \nonumber\\
&&\ds=\int_{\t}^{\t+\theta}\int_{\t}^{t}
\dbE~\Big\{\Big\langle e^{A(t-s)}\[ b_1(s)y_{ij}^{k,\th}(s)+ b_2(s)\big(v^{k}-\bar{u}(s)\big)\chi_{F_{ij}}(\omega)\], \nonumber\\
&&\ds
\qquad\qquad\qquad\qquad\qquad\qquad\qquad\qquad
\cD_{s}\big(\dbS(t)^{*}(v^{k}-\bar{u}(t))\big)
\Big\rangle_{H}\chi_{F_{ij}}(\omega)\Big\}dsdt.\nonumber
\end{eqnarray}

By Lemma \ref{lm4}, there exists a sequence
$\{\theta_{n}\}_{n=1}^{\infty}$ such that
$\theta_n\to0^+$ as $n\to\infty$ and
\begin{equation}\label{12.20-eq12}
\begin{array}{ll}\ds
\q\frac{1}{\theta_{n}^2}\Big|\int_{\t}^{\t+\theta_{n}}\int_{\t}^{t}
\dbE~\Big[\Big\langle e^{A(t-s)}
b_1(s)y_{ij}^{k,\th_{n}}(s),
\cD_{s}\big(\dbS(t)^{*}(v^{k}-\bar{u}(t))\big)
\Big\rangle_{H}\chi_{F_{ij}}(\omega)\Big]dsdt\Big|\\
\ns\ds \le \frac{1}{\theta_{n}^2}\Big(\dbE\!
\int_{\t}^{\t+\theta_{n}}\!\!\!\!
\int_{\t}^{t}\big|e^{A(t-s)}
b_1(s)y_{ij}^{k,\th_{n}}(s)
\big|_{H}^2dsdt\Big)^{\frac{1}{2}}\!\cdot
\Big(\dbE\int_{\t}^{\t+\theta_{n}}\!\!\!\int_{\t}^{t}
\big|\cD_{s}\big(\dbS(t)^{*}(v^{k}\!\!-\!\bar{u}(t))\big)
\big|^2_Hdsdt\Big)^{\frac{1}{2}}\\
\ns\ds \le \frac{C}{\theta_{n}}
\Big(\sup_{t\in[0,T]}\dbE|y_{ij}^{k,\th_{n}}(t)|_{H}^{2}\Big)^{\frac{1}{2}}\cdot
\Big(\dbE~
\int_{\t}^{\t+\theta_{n}}\!\int_{\t}^{t}\!
\big|\cD_{s}\big(\dbS(t)^{*}(v^{k}\!-\!\bar{u}(t))\big)
\big|^2_Hdsdt\Big)^{\frac{1}{2}}\\
\ns\ds\q\to 0,\q n\to\infty, \q a.e.\ \t\in[t_{i},T).
\end{array}
\end{equation}

We next prove that there exists a subsequence
$\{\theta_{n_l}\}_{l=1}^{\infty}$ of
$\{\theta_{n}\}_{n=1}^{\infty}$,  such that
$\theta_{n_l}\to 0^+$ as $l\to\infty$ and
\begin{equation}\label{12.20-eq13}
\begin{array}{ll}\ds
\q\lim_{l\to \infty}\frac{1}{\theta_{n_l}^{2}}
\int_{\t}^{\t+\theta_{n_l}}\int_{\t}^{t}
\dbE~\Big(\big\langle e^{A(t-s)}b_2(s)\big(v^{k}-\bar{u}(s)\big)\chi_{F_{ij}}(\omega),\\
\ns\ds
\qquad\qquad\qquad\qquad\qquad\qquad\qquad\qquad
\cD_{s}\big(\dbS(t)^{*}(v^{k}-\bar{u}(t))\big)
\big\rangle_{H}\chi_{F_{ij}}(\omega)\Big)dsdt.\\
\ns\ds=\frac{1}{2} \dbE~\big(\big\langle
b_{2}(\t) (v^{k}-\bar{u}(\t)),\nabla
\dbS(\t)^*\big(v^{k}-\bar{u}(\t)\big)\big\rangle_{H}\chi_{F_{ij}}(\omega)\big) \\
\ns\ds\q -\frac{1}{2} \dbE~\big(\langle
b_{2}(\t) (v^{k}-\bar{u}(\t)),
\dbS(\t)^*\nabla\bar{u}(\t)\rangle_{H}\chi_{F_{ij}}(\omega)\big),\
\ a.e. \ \t\in [t_{i},T).
\end{array}
\end{equation}

By \textbf{(A6)},
$$
\cD_{s}\big(\dbS(t)^{*}(v^{k}-\bar{u}(t))\big)
=\cD_{s} \dbS(t)^{*}(v^{k}-\bar{u}(t))
-\dbS(t)^{*}\cD_{s}\bar{u}(t).
$$
Then, we have
\begin{equation}\label{12.20-eq14}
\begin{array}{ll}\ds
\q\frac{1}{\theta^2}\int_{\t}^{\t+\theta}\int_{\t}^{t}
\dbE~\Big(\big\langle e^{A(t-s)}b_{2}(s)
(v^{k}-\bar{u}(s))\chi_{F_{ij}}(\omega), \\
\ns\ds
\qquad\qquad\qquad\qquad\qquad\qquad\qquad\qquad\
\ \
\cD_{s}\big(\dbS(t)^{*}(v^{k}-\bar{u}(t))\big)
\big\rangle_{H}\chi_{F_{ij}}(\omega)\Big)dsdt \\
\ns\ds=\frac{1}{\theta^2}\int_{\t}^{\t+\theta}\int_{\t}^{t}
\dbE~\Big(\big\langle e^{A(t-s)}b_{2}(s)
(v^{k}-\bar{u}(s)),  \\
\ns\ds
\qquad\qquad\qquad\qquad\qquad\qquad\qquad\qquad\
\ \ \ \ \cD_{s} \dbS(t)^{*}(v^{k}-\bar{u}(t))
\big\rangle_{H}\chi_{F_{ij}}(\omega)\Big)dsdt \\
\ns\ds\q
-\frac{1}{\theta^2}\int_{\t}^{\t+\theta}\int_{\t}^{t}
\dbE~\Big(\big\langle e^{A(t-s)}b_{2}(s)
(v^{k}-\bar{u}(s)),\dbS(t)^{*}\cD_{s}\bar{u}(t)
\big\rangle_{H}\chi_{F_{ij}}(\omega)\Big)dsdt.
\end{array}
\end{equation}

For the first part in the right hand side of
(\ref{12.20-eq14}),
\begin{equation}\label{12.20-eq15}
\begin{array}{ll}\ds
\q\frac{1}{\theta^2}\int_{\t}^{\t+\theta}\int_{\t}^{t}
\dbE~\Big(\big\langle e^{A(t-s)}b_{2}(s)
(v^{k}-\bar{u}(s)),  \cD_{s}
\dbS(t)^{*}(v^{k}-\bar{u}(t))\big\rangle_{H}
\chi_{F_{ij}}(\omega)\Big)dsdt \\
\ns\ds=\frac{1}{\theta^2}\int_{\t}^{\t+\theta}\!\!\int_{\t}^{t}\!
\dbE~\Big[\big\langle e^{A(t-s)}b_{2}(s)
(v^{k}\!-\!\bar{u}(s)),  \big(\cD_{s}
\dbS(t)\!-\!\nabla\dbS(s)\big)^{*}
(v^{k}\!-\!\bar{u}(t))\big\rangle_{H}\chi_{F_{ij}}(\omega)\Big]dsdt \\
\ns\ds\q
+\frac{1}{\theta^2}\int_{\t}^{\t+\theta}\int_{\t}^{t}
\dbE~\Big(\big\langle e^{A(t-s)}b_{2}(s)
(v^{k}-\bar{u}(s)),
\nabla\dbS(s)^{*}(v^{k}-\bar{u}(t))\big\rangle_{H}
\chi_{F_{ij}}(\omega)\Big)dsdt.
\end{array}
\end{equation}
Since
$$
\begin{array}{ll}\ds
\Big|\frac{1}{\theta^2}\int_{\t}^{\t+\theta}\int_{\t}^{t}
\dbE\Big[\big\langle e^{A(t-s)}b_{2}(s)
(v^{k}-\bar{u}(s)),\\
\ns\ds \qquad\qquad\qquad\qquad\qquad
\big(\cD_{s} \dbS(t)-\nabla\dbS(s)\big)^{*}
(v^{k}-\bar{u}(t))\big\rangle_{H}\chi_{F_{ij}}(\omega)\Big]dsdt\Big|\\
\ns\ds\le
\frac{C}{\theta^2}\int_{\t}^{\t+\theta}\int_{\t}^{t}\dbE~
\Big[\big|e^{A(t-s)}b_{2}(s)
(v^{k}-\bar{u}(s))\big|_{H}\cdot\\
\ns\ds\qq\qq\qq\qq\qq\qq\qq
\big|\cD_{s} \dbS(t)-\nabla\dbS(s)\big|_{\cL_2(H_1;H)}\cdot|v^{k}-\bar u(t)|_{H_{1}}\Big]dsdt\\
\ns\ds\le\!\frac{C}{\theta^2}
\(\dbE\!\int_{\t}^{\t+\theta}\!\!\int_{\t}^{t}\!\!
\big|
(v^{k}\!-\!\bar{u}(s))\big|_{H_{1}}^2\cdot\big|
(v^{k}\!-\!\bar{u}(t))\big|_{H_{1}}^2 dsdt\)^{\frac{1}{2}}\!\cdot\!\\
\ns\ds\qq\qq\qq\qq\qq\qq\qq\q
\(\dbE\!\int_{\t}^{\t+\theta}\!\int_{\t}^{t}\!\!
\big|\cD_{s} \dbS(t)\!-\!\nabla\dbS(s)\big|_{\cL_2(H_1;H)}^2dsdt\)^{\frac{1}{2}},
\end{array}
$$
by Lemma \ref{lm4}, there exists a a subsequence
of $\{\theta_{n_l^1}\}_{l=1}^{\infty}$ of
$\{\theta_{n}\}_{n=1}^{\infty}$ such that
$\theta_{n_l^1}\to0^+$ as $l\to\infty$ and
\begin{equation}\label{12.20-eq16}
\begin{array}{ll}\ds
\lim_{l\to
\infty}\frac{1}{\theta_{n_l^1}^{2}}\int_{\t}^{\t+\theta_{n_l^1}}\int_{\t}^{t}
\dbE\Big[\big\langle e^{A(t-s)}b_{2}(s)
(v^{k}-\bar{u}(s)),\\
\ns\ds \qquad\qquad\qquad\qquad\qquad
\big(\cD_{s}\dbS(t)-\nabla\dbS(s)\big)
^{*}(v^{k}-\bar{u}(t))\big\rangle_{H}\chi_{F_{ij}}(\omega)\Big]dsdt \\
\ns\ds =0,\qq \qquad a.e. \;\t\in[0,T).
\end{array}
\end{equation}

For the second part in the right hand side of
(\ref{12.20-eq15}), by Lemma \ref{lm3} it
follows that
\begin{equation}\label{12.20-eq17}
\begin{array}{ll}\ds
 \lim_{\theta\to
0^+}\frac{1}{\theta^2}\int_{\t}^{\t+\theta}\int_{\t}^{t}
\dbE~\Big(\big\langle e^{A(t-s)}b_{2}(s)
(v^{k}-\bar{u}(s)),
\nabla\dbS(s)^{*}(v^{k}-\bar{u}(t))\big\rangle_{H}
\chi_{F_{ij}}(\omega)\Big)dsdt \\
\ns\ds =\frac{1}{2}\dbE~\big(\big\langle
b_{2}(\t)(v^{k}
-\bar{u}(\t)),\nabla\dbS(\t)^*\big(v^{k}-\bar{u}(\t)\big)\big\rangle_{H}\chi_{F_{ij}}(\omega)
\big),\qq \qq a.e. \ \t\in [t_{i},T).
\end{array}
\end{equation}

Therefore, by
(\ref{12.20-eq15})--(\ref{12.20-eq17}), we
conclude that
\begin{equation}\label{12.20-eq18}
\begin{array}{ll}\ds
\q\lim_{l\to
\infty}\frac{1}{\theta_{n_l^1}^{2}}\int_{\t}^{\t+\theta_{n_l^1}}\int_{\t}^{t}
\dbE~\Big(\big\langle e^{A(t-s)}b_{2}(s)
(v^{k}-\bar{u}(s)), \cD_{s}
\dbS(t)^{*}(v^{k}-\bar{u}(t))\big\rangle_{H}
\chi_{F_{ij}}(\omega) \Big)dsdt \\
\ns\ds=\frac{1}{2}\dbE~\big(\langle
b_{2}(\t)(v^{k}
-\bar{u}(\t)),\nabla\dbS(\t)^*\big(v^{k}-\bar{u}(\t)\big)\rangle_{H}\chi_{F_{ij}}(\omega)\big),\
\qq\qq a.e. \ \t\in [t_{i},T).
\end{array}
\end{equation}

In a similar way, we can prove that there exists
a subsequence
$\{\theta_{n_{l}}\}_{l=1}^{\infty}$ of
$\{\theta_{n_l^1}\}_{n=1}^{\infty}$  such that
\begin{equation}\label{12.20-eq19}
\begin{array}{ll}\ds
\q\lim_{l\to
\infty}\frac{1}{\theta_{n_{l}}^{2}}\int_{\t}^{\t+\theta_{n_{l}}}\int_{\t}^{t}
\dbE~\Big(\big\langle e^{A(t-s)}b_{2}(s)
(v^{k}-\bar{u}(s)),\dbS(t)^{*}\cD_{s}\bar{u}(t)
\big\rangle_{H}\chi_{F_{ij}}(\omega) \Big)dsdt \\
\ns\ds= \frac{1}{2} \dbE~\big(\langle
b_{2}(\t)(v^{k}-\bar{u}(\t)),
\dbS(\t)^*\nabla\bar{u}(\t)\rangle_{H}\chi_{F_{ij}}(\omega)
\big),\qq \ a.e. \ \t\in [t_{i},T).
\end{array}
\end{equation}

Combining (\ref{12.20-eq14}), (\ref{12.20-eq18})
and (\ref{12.20-eq19}),  we obtain
(\ref{12.20-eq13}). Then, by (\ref{12.20-eq11}),
(\ref{12.20-eq12}) and (\ref{12.20-eq13}),  we
obtain that there exists a subsequence
$\{\theta_{n_l}\}_{l=1}^{\infty}$,
$\theta_{n_l}\to 0^+$ as $l\to \infty$ and
\begin{equation}\label{12.20-limit part2}
\begin{array}{ll}\ds
\q\lim_{l\to
\infty}\dbE\int_{\t}^{\t+\theta_{n_l}}
\Big\langle \int_{\t}^{t}
e^{A(t-s)}\[ b_1(s)y_{ij}^{k,\th_{n_l}}(s)\\
\ns\ds\qquad\qquad\qquad\qquad\qquad
+ b_2(s)\big(v^{k}-\bar{u}(s)\big)\chi_{F_{ij}}(\omega)
\]dW(s),
\dbS(t)\big(v^{k}-\bar{u}(t)\big)\Big\rangle_{H}\chi_{F_{ij}}(\omega) dt\\
\ns\ds=\frac{1}{2} \dbE~\big(\langle b_{2}(\t)
(v^{k}-\bar{u}(\t)),\nabla
\dbS(\t)^*\big(v^{k}-\bar{u}(\t)\big)\rangle_{H}\chi_{F_{ij}}(\omega) \big) \\
\ns\ds\q -\frac{1}{2} \dbE~\big(\langle
b_{2}(\t) (v^{k}-\bar{u}(\t)),
\dbS(\t)^*\nabla\bar{u}(\t)\rangle_{H}\chi_{F_{ij}}(\omega)
\big),\ \ a.e. \ \t\in [t_{i},T).
\end{array}
\end{equation}

Finally, by (\ref{12.20-eq9}),
(\ref{12.20-limit part1}) and (\ref{12.20-limit part2}) we
conclude that, for any $i,j,k\in \dbN$, there
exists a Lebesgue measurable set
$E^{k}_{i,j}\subset[t_{i},T)$ with
$|E^{k}_{i,j}|=0$ such that
\begin{equation}\label{12.20-eq27}
\begin{array}{ll}\ds
 0\3n&\ds\ge \frac{1}{2}\dbE~\big(
\langle a_{2}(\t)(v^{k}
-\bar{u}(\t)),\dbS(\t)^*\big(v^{k}-\bar{u}(\t)\big)\rangle_{H}\chi_{F_{ij}}(\omega)\big)\\
\ns&\ds \q+\frac{1}{2} \dbE~\big(\langle
b_{2}(\t)(v^{k}-\bar{u}(\t)),\nabla
\dbS(\t)^*\big(v^{k}-\bar{u}(\t)\big)
\rangle_{H}\chi_{F_{ij}}(\omega)\big)\nonumber\\
\ns&\ds \q -\frac{1}{2} \dbE~\big(\langle
b_{2}(\t)(v^{k}-\bar{u}(\t)),\dbS(\t)^*\nabla\bar{u}(\t)
\rangle_{H}\chi_{F_{ij}}(\omega)\big),\quad
\forall \ \t\in [t_{i},T)\setminus E^{k}_{i,j}.
\end{array}
\end{equation}
Let $E_{0}=\bigcup_{i,j,k\in \dbN} E^{k}_{i,j}$.
Then $|E_{0}|=0$, and for any $i,j,k\in \dbN$,
\begin{equation*}
\begin{array}{ll}\ds
 \dbE\big( \langle a_{2}(\t)(v^{k}
-\bar{u}(\t)),\dbS(\t)^*\big(v^{k}-\bar{u}(\t)\big)\rangle_{H}\chi_{F_{ij}}(\omega)\big)\\[+0.4em]
\ns\ds\quad+ \dbE~\big( \langle
b_{2}(\t)(v^{k}-\bar{u}(\t)),\nabla
\dbS(\t)^*\big(v^{k}-\bar{u}(\t)\big)
\rangle_{H}\chi_{F_{ij}}(\omega)\big) \\[+0.4em]
\ns\ds\quad-\dbE~\big( \langle
b_{2}(\t)(v^{k}-\bar{u}(\t)),\dbS(\t)^*\nabla\bar{u}(\t)
\rangle_{H}\chi_{F_{ij}}(\omega)\big) \\[+0.4em]
\ns\ds\le 0, \  \  \forall\ \t\in
[t_{i},T)\setminus E_{0}.
\end{array}
\end{equation*}
By the construction of
$\{F_{ij}\}_{i=1}^{\infty}$, the continuity of
the filter $\dbF$ and the density of
$\{v^{k}\}_{k=1}^{\infty}$, we conclude that
\begin{equation*}
\begin{array}{ll}\ds
\langle a_{2}(\t)(v
-\bar{u}(\t)),\dbS(\t)^*\big(v-\bar{u}(\t)\big)\rangle_{H}
+ \langle
 b_{2}(\t)(v-\bar{u}(\t)),\nabla\dbS(\t)^*\big(v-\bar{u}(\t)\big)\rangle_{H} \\[+0.4em]
\ns\ds -\langle
b_{2}(\t)(v-\bar{u}(\t)),\dbS(\t)^*\nabla\bar{u}(\t)\rangle_{H}
 \le 0,\ \ a.s., \qquad   \forall\ (\t,v)\in
([0,T]\setminus E_{0})\times U.
\end{array}
\end{equation*}
This completes the proof of Theorem \ref{th max
p1}.
\endpf

%%%%%%%%%%%%%%%%%%%%%%%%%%%%%%%%%%%%%%%%%%%%

\section{Second order sufficient
conditions}\label{s5}

In this section, we discuss the second order
sufficient condition for the optimal control
problem (\ref{jk2}). We first give a simple and
direct result, and then we generalize it under
some proper assumptions and obtain a second
order sufficient condition which has minimal gap
with the second order necessary condition. The
basic idea comes from  optimization theory.

In addition to Assumption
\textbf{(A1)}--\textbf{(A3)}, we assume that \ms

\no{\bf (A7)} {\it  $U$ is a bounded  closed convex set.}

\no{\bf (A8)} {\it There exists a constant
$C_L>0$ such that, for any $t\in[0,T]$ and
$(x,u), (\tilde x, \tilde u) \in H\times U$,
\begin{equation}\label{Lipschitz seccond order deriv}
\left\{
\begin{array}{ll}\ds
\|a_{xx}(t,x,u)-a_{xx}(t,\tilde x,\tilde
u)\|_{\cL(H\times
H;H)}+\|b_{xx}(t,x,u)-b_{xx}(t,\tilde x,\tilde
u)\|_{\cL(H\times
H;H)}\\
\ns\ds +\|a_{xu}(t,x,u)-a_{xu}(t,\tilde x,\tilde
u)\|_{\cL(H\times
H_{1};H)}+\|b_{xu}(t,x,u)-b_{xu}(t,\tilde
x,\tilde u)\|_{\cL(H\times
H_{1};H)}\\
\ns\ds +\|a_{uu}(t,x,u)-a_{uu}(t,\tilde x,\tilde u)\|_{\cL(H_{1}\times H_{1};H)}+\|b_{uu}(t,x,u)-b_{uu}(t,\tilde x,\tilde u)\|_{\cL(H_{1}\times H_{1};H)}\\
\ns\ds
\leq C_L\big(|x-\tilde x|_{H} + |u-\tilde u|_{H_{1}}\big),\\
\ns\ds
\|g_{xx}(t,x,u)-g_{xx}(t,\tilde x,\tilde u)\|_{\cL(H)}
+\|g_{xu}(t,x,u)-g_{xu}(t,\tilde x,\tilde u)\|_{\cL(H;H_{1})}\\
\ns\ds
+\|g_{uu}(t,x,u)-g_{uu}(t,\tilde x,\tilde u)\|_{\cL(H_{1})}
+\|h_{xx}(x)-h_{xx}(\tilde x)\|_{\cL(H)}\\
\ns\ds \leq C_L\big(|x-\tilde x|_{H} + |u-\tilde u|_{H_{1}}\big).
\end{array}
\right.
\end{equation}}

Under Assumption \textbf{(A7)}, any $U$-valued
measurable adapted process $u(\cdot)$ belongs to
$\cU^{\infty}[0,T]\subset\cU^{\beta}[0,T]$
($\beta\ge2$). Let $u(\cdot), \bar
u(\cdot)\in\cU^{\infty}[0,T]$, $x(\cdot)$ and
$\bar x(\cdot)$  be solutions to the control
system  (\ref{fsystem1}) with respect to
$u(\cdot)$ and $\bar u(\cdot)$, respectively.
Let $\d u$, $\d x$, $y$ and $z$ be defined as in
Section 2. We first give the following estimate:
\begin{lemma}\label{lemma 5.2}
Let {\bf(A1)},  {\bf(A3)} and
{\bf(A7)}--{\bf(A8)} hold. Then, for any $\b\ge
2$,
\begin{equation}\label{th suff estmate r2}
\Big\|\delta
x-y-\frac{1}{2}z\Big\|_{C_\dbF([0,T];L^{\b}(\O;H))}
\le C\big(\|\d
u\|_{L^{\infty}_\dbF(0,T;H_1)}\cdot\|\d
u\|_{L^{4\b}_\dbF(0,T;H_1)}^2\big).
\end{equation}
\end{lemma}
{\it Proof}\,: For $\psi=a,b$, put
$$
\begin{cases}\ds
\hat{\psi}_{11}(t)\=\int_{0}^{1}(1-\theta)\psi_{xx}(t,\bar
x(t)+\theta\delta x(t),\bar u(t)+\theta\d
u(t))d\theta,\\
\ns\ds
\hat{\psi}_{12}(t)\=\int_{0}^{1}(1-\theta)\psi_{xu}(t,\bar
x(t)+\theta\delta x(t),\bar u(t)+\theta\d
u(t))d\theta, \\
\ns\ds
\hat{\psi}_{22}(t)\=\int_{0}^{1}(1-\theta)\psi_{uu}(t,\bar
x(t)+\theta\delta x(t),\bar u(t)+\theta\d
u(t))d\theta.
\end{cases}
$$

Similar to the proof of Step 1 in Theorem
\ref{th max}, $\d x$ solves the following SEE:
\begin{equation}\label{th suff deltax order two}
\left\{
\begin{array}{ll}
d\delta x= \Big[A\delta x +a_1(t)\delta x
+a_{2}(t)\delta u
+\hat{a}_{11}(t)\big(\delta x,\delta x\big)\\[+0.4em]
\qquad\qquad\qq\qq\qq\ \
+2\hat{a}_{12}(t)\big(\delta x,\delta u\big)
+\hat{a}_{22}(t)\big(\delta u,\delta u\big)
\Big]dt\\[+0.4em]
\qquad\qquad
+\Big[b_1(t)\delta x
+ b_{2}(t)\delta u
+\hat{b}_{11}(t)\big(\delta x,\delta x\big)\\[+0.4em]
\qquad\qquad\qquad\
+2\hat{b}_{12}(t)\big(\delta x,\delta u\big)
+ \hat{b}_{22}(t)\big(\delta u,\delta u\big)\Big]dW(t) &\mbox{ in } (0,T],\\[+0.4em]
\delta x(0)=0.
\end{array}\right.
\end{equation}

Let $r_{2}(\cd) =\d
x(\cdot)-y(\cdot)-\frac{1}{2}z(\cdot)$. Then
$r_{2}(\cd)$ fulfills
\begin{equation}\label{th suff-eq3}
\left\{
\begin{array}{ll}\ds
dr_{2}=\big(Ar_{2} + a_1(t)r_{2} +
\Upsilon_{1}(t)\big)dt   + \big(b_1(t)r_{2} +
\Upsilon_{2}(t) \big)dW(t) \ \ \mbox{ in
}(0,T],\\
\ns\ds r_{2}(0)=0,
\end{array}
\right.
\end{equation}
where
$$
\begin{array}{ll}\ds
\Upsilon_{1}(t)\ds=\( \hat a_{11}(t)(\d x(t),\d
x(t))-
\frac{1}{2}a_{11}(t)(y(t),y(t))\)\\
\ns\ds \qq\qq  + \big(2\hat a_{12}(t)(\d x(t),\d
u(t))\! -\! a_{12}(t)(y(t),\d u(t))\big) \!+\!
\(\hat a_{22}(t)\!-\!\frac{1}{2}a_{22}(t)\)(\d
u(t),\d u(t))
\end{array}
$$
and
$$
\begin{array}{ll}\ds
\Upsilon_{2}(t)\ds=\(\hat b_{11}(t)(\d x(t),\d
x(t)) -\frac{1}{2} b_{11}(t)(y(t),y(t))\)
\\
\ns\ds \qq\qq +  \big(2\hat b_{12}(t)(\d x(t),\d
u(t)) \!-\! b_{12}(t)(y(t),\d u(t))\big) \!+\!
\(\hat b_{22}(t)\!-\!\frac{1}{2}b_{22}(t)\)(\d
u(t),\d u(t)).
\end{array}
$$
We have that
\begin{equation}\label{th suff-eq7}
\begin{array}{ll}\ds
\mE\big|r_2(t)\big|_H^{\b}= \mE\Big| \int_0^t
e^{A(t-s)} a_1(s)r_2(s)ds
+ \int_0^t e^{A(t-s)} b_1(s)r_2(s) dW(s) \\
\ns\ds\qq\qq\qq +\int_0^t
e^{A(t-s)}\Upsilon_{1}(s) ds +\int_0^t
e^{A(t-s)}\Upsilon_{2}(s) dW(s) \Big|_H^{\b}
\\
\ns\ds\ \qq\qq \leq C\[\mE\int_0^t
|r_2(s)|_H^{\b} ds + \mE\(\int_0^t
|\Upsilon_{1}(s)|_H ds\)^{\b} + \mE\(\int_0^t
|\Upsilon_{2}(s)|_H^{2} ds\)^{\frac{\b}{2}}\].
\end{array}
\end{equation}
By \textbf{(A3)}, \textbf{(A7)}--\textbf{(A8)}
and Lemma \ref{main estmates}, we deduce that
\begin{eqnarray}\label{th suff-eq8}
&&\mE\(\int_0^t
|\Upsilon_{1}(s)|_Hds\)^{\b}\nonumber\\
&\le&\3n C\mE\(\int_0^T\Big| \big( \hat
a_{11}(t)(\d x(t),\d x(t))-
\frac{1}{2}a_{11}(t)(y(t),y(t))\big)\nonumber\\
&& \qq\qq   + \big(2\hat a_{12}(t)(\d x(t),\d
u(t)) -
a_{12}(t)(y(t),\d u(t))\big) \nonumber\\
&& \qq\qq
+\big(\hat a_{22}(t)-\frac{1}{2}a_{22}(t)\big)(\d u(t),\d u(t))\Big|_H dt\)^{\b}\nonumber\\
&\le&\3n\! C\mE\[\!\int_0^T\!\!\(\big| \hat
a_{11}(t)(\d x(t),\d x(t))\!-\! \hat
a_{11}(t)(y(t),y(t))|_{H}\! +\!\big\|\hat
a_{11}(t)\!-\!
\frac{1}{2}a_{11}(t)\|_{\cL(H\times H;H)} \cdot|y(t)|_H^{2}\nonumber\\
&& \qq\qq  +2\big|\hat a_{12}(t)(\d x(t),\d
u(t)) -
\hat a_{12}(t)(y(t),\d u(t))\big|_H \nonumber\\
&&\qq\qq  +\big\|2\hat a_{12}(t) -
a_{12}(t)\big\|_{\cL(H\times H_{1};H)}
\cdot|y(t)|_H \cdot|\d
u(t)|_{H_{1}} \nonumber \\
&&\qq\qq  +
\big\|\hat a_{22}(t)-\frac{1}{2}a_{22}(t)\big\|_{\cL(H_{1}\times H_{1};H)} \cdot|\d u(t)|_{H_{1}}^{2}\)dt\]^{\b}\nonumber\\
&\leq&\3n\! C\mE\[\int_0^T\(\big|\d
x(t)+y(t)\big|_{H} \cdot\big|\d x(t)-y(t)|_{H}
+\big(|\d x(t)|_{H} +|\d u(t)|_{H_{1}}  \big)\cdot|y(t)|_{H}^{2}\nonumber\\
&& \qq\qq  +2\big|\d x(t) -y(t)\big|_{H}
\cdot\big|\d u(t)\big|_{H_{1}}  +\big(|\d
x(t)|_{H} +|\d
u(t)|_{H_{1}} \big)\cdot|y(t)|_{H} \cdot|\d u(t)|_{H_{1}} \nonumber\\
&& \qq\qq  +
\big(|\d x(t)|_{H} +|\d u(t)|_{H_{1}} \big)\cdot|\d u(t)|_{H_{1}}^{2}\)dt\]^{\b}\nonumber\\
& \le&\3n C \[\mE\(\int_0^T|\d
u(t)|_{H_{1}}^{2}dt\)^{\b}
\]^{\frac{1}{2}} \cd \[\mE\(\int_0^T|\d
u(t)|_{H_{1}}^{4}dt\)^{\b}
\]^{\frac{1}{2}}\nonumber\\
&\le&\3n\! C \|\d
u\|_{L^{\infty}_\dbF(0,T;H_1))}^{\b}\cdot\|\d
u\|_{L^{4\b}_\dbF(0,T;H_1)}^{2\b}.
\end{eqnarray}
Similarly,
\begin{equation}\label{th suff-eq9}
\mE\(\int_0^t |\Upsilon_{2}(s)|_H^{2}
ds\)^{\frac{\b}{2}}\le  C \|\d u\|_{L^{\infty}_\dbF(0,T;H_1))}^{\b}\cdot\|\d
u\|_{L^{4\b}_\dbF(0,T;H_1)}^{2\b}.
\end{equation}
Combining \eqref{th suff-eq7}, \eqref{th suff-eq8} with \eqref{th suff-eq9}, we obtain \eqref{th suff estmate r2}.
\endpf

\vspace{+0.5em}

Now, we put
\begin{equation}\label{Lambda}
\begin{array}{ll}\ds
\Lambda (v(\cdot))\=\mE\int_0^T \[
\big\langle \dbH_{uu}(t)v(t), v(t) \big\rangle_{H_{1}}\! +\!
\big\langle
b_2(t)^{*}P(t)b_2(t)v(t), v(t)\big\rangle_{H_{1}}
\]dt\\
\ns\ds\qq\qq\ +2\mE\!\int_0^T\!\!\!
\big\langle\big(\dbH_{xu}(t)\! +\! a_2(t)^*
P(t)\! +\! b_2(t)^*P(t) b_1(t)\big)y^{v}(t), v(t)\big\rangle_{H} dt  \\
\ns\ds\qq\qq\ + \mE\int_0^T
\big\langle\big( \widehat Q^{(0)}+Q^{(0)}\big)\big(0,a_2(t)v(t),
b_2(t)v(t)\big), b_2(t)v(t)
\big\rangle_{H}dt,
\end{array}
\end{equation}
and
\begin{equation}\label{tildeLambda}
\begin{array}{ll}\ds
\tilde\Lambda (v(\cdot))\=\mE\int_0^T \[
\big\langle \dbH_{uu}(t)v(t), v(t) \big\rangle_{H_{1}}\! +\!
\big\langle
b_2(t)^{*}P(t)b_2(t)v(t), v(t)\big\rangle_{H_{1}}
\]dt\\
\ns\ds\qq\qq\ +2\mE\!\int_0^T\!\!\!
\big\langle\big(\dbS(t)y^{v}(t), v(t)\big\rangle_{H_{1}} dt,
\end{array}
\end{equation}
where $y^{v}(\cdot)$ is the solution to the
equation (\ref{fsystem3.1}) with $\d u$ replaced
by $v$ and $(P(\cd),(Q^{(\cd)},\widehat
Q^{(\cd)}))$ (resp. $(P,Q)$) is the relaxed
transposition solution (resp. the
$V$-transposition solution) of BSEE
\eqref{op-bsystem3} with $P_T$, $J(\cd)$,
$K(\cd)$ and $F(\cd)$ given by (\ref{zv2}). Note
that the mapping
$$
v(\cdot)\mapsto y^{v}(\cdot)
$$
from ${L^{\b}_\dbF(0,T;H_1)}$ to
$C_{\dbF}([0,T];L^{\beta}(\O;H))$ is a linear
for any $\b\ge 2$. The mapping $\Lambda$ and
$\tilde\Lambda$ are actually two quadratic-like
forms defined on the Banach space
${L^{\b}_\dbF(0,T;H_1)}$ for any $\b\ge 4$.

By Lemmas \ref{main estmates} and \ref{lemma 5.2}, we obtain the following second
order sufficient condition.
\begin{theorem}\label{th second order suff}
Suppose that $x_0\in L^2_{\cF_0}(\O;H)$. Let
{\bf(A1)}--{\bf(A3)} and {\bf(A7)}--{\bf(A8)}
hold, and let $\bar u(\cd)$ be an admissible
control and $\bar x(\cd)$ be the corresponding
state. If there exists a constant $\varrho>0$
such that for any $u(\cd) \in
\cU^{\infty}[0,T]$,
\begin{equation}\label{1st order condition}
\mE\int_{0}^{T}\big\langle\dbH_{u}(t), u(t)-\bar
u(t)  \big\rangle_{H_{1}} dt\le0
\end{equation}
and
\begin{equation}\label{2nd order condition}
\Lambda(u-\bar u) \leq - 2\varrho\|u-\bar u\|_{L^{8}_\dbF(0,T;H_1)}^2,
\end{equation}
then there exists a constant $\sigma>0$ such that for any $u(\cd) \in \cU^{\infty}[0,T]$ with $\|u-\bar u\|_{L^{\infty}_\dbF(0,T;H_1)}\le \sigma$,
\begin{equation}\label{strong local min}
\cJ(u)\ge \cJ(\bar u)+\frac{\varrho}{2}\|u-\bar u\|_{L^{8}_\dbF(0,T;H_1)}^2.
\end{equation}
Especially, $\bar u$ is a local  minima of the
optimal control problem (\ref{jk1}).
\end{theorem}
{\it Proof}\,: Let $u(\cdot)\in
\cU^{\infty}[0,T]$ and $x(\cdot)$(\resp $\bar
x(\cdot)$) be the solutions to the control
system (\ref{fsystem1}) with respect to
$u(\cdot)$(\resp $\bar u(\cdot)$). Let $\d u$,
$\d x$, $y$ and $z$ be defined as in Section
\ref{s2}. Put
$$
\begin{cases}\ds
\hat{g}_{11}(t)\=\int_{0}^{1}(1-\theta)g_{xx}(t,\bar
x(t)+\theta\delta x(t),\bar u(t)+\theta \d
u(t))d\theta,\\
\ns\ds
\hat{g}_{12}(t)\=\int_{0}^{1}(1-\theta)g_{xu}(t,\bar
x(t)+\theta\delta x(t),\bar u(t)+\theta \d
u(t))d\theta,\\
\ns\ds
\hat{g}_{22}(t)\=\int_{0}^{1}(1-\theta)g_{uu}(t,\bar
x(t)+\theta\delta x(t),\bar u(t)+\theta \d
u(t))d\theta,\\
\ns\ds \hat{h}_{xx}(T)
\=\int_{0}^{1}(1-\theta)h_{xx}(\bar
x(T)+\theta\delta x(T))d\theta.
\end{cases}
$$
By Taylor's formula, we see that
\begin{equation}\label{th suff-eq10}
\begin{array}{ll}\ds
\q g(t,x(t),u(t)) - g(t,\bar x(t),
\bar u(t)) \\[+0.4em]
\ns\ds = \big\langle g_1(t),  \delta x(t)
\big\rangle_H + \big\langle g_2(t),  \d u(t)
\big\rangle_{H_{1}} + \big\langle
\hat{g}_{11}(t)\delta x(t), \delta x(t) \big\rangle_H \\[+0.4em]
\ns\ds \q +2\big\langle \hat{g}_{12}(t)\delta
x(t), \d u(t) \big\rangle_{H_{1}} +\big\langle
\hat{g}_{22}(t)\d u(t), \d u(t)
\big\rangle_{H_{1}},
\end{array}
\end{equation}
and
\begin{equation}\label{th suff-eq11}
 h(x(T)) - h(\bar x(T))
= \big\langle h_x(\bar x(T)), \d x(T)
\big\rangle_H + \big\langle
\hat{h}_{xx}(T)\d x(T), \d x(T)
\big\rangle_H.
\end{equation}

Using a similar method in the proof  of (\ref{th
suff-eq8}), we obtain that
\begin{equation}\label{th suff-eq12}
\begin{array}{ll}\ds
\q \Big|
\dbE \int_{0}^{T}\Big(\big\langle\hat{g}_{11}(t)\delta x(t),\delta x(t)\big\rangle_H-\frac{1}{2}\big\langle g_{11}(t)y(t), y(t)\big\rangle_H\Big)dt\Big|\\[+0.4em]
\ns\ds\le C\big(\|\d
u\|_{L^{\infty}_\dbF(0,T;H_1)}\cdot\|\d
u\|_{L^{4}_\dbF(0,T;H_1)}^2\big)\le  C\big(\|\d
u\|_{L^{\infty}_\dbF(0,T;H_1)}\cdot\|\d
u\|_{L^{8}_\dbF(0,T;H_1)}^2\big) ,
\end{array}
\end{equation}
\begin{equation}\label{th suff-eq13}
\begin{array}{ll}\ds
\q \Big|
\dbE \int_{0}^{T}\Big(2\big\langle\hat{g}_{12}(t)\delta x(t), \d u(t)\big\rangle_{H_{1}}- \big\langle g_{12}(t)y(t), \d u(t)\big\rangle_{H_{1}}\Big)dt\Big|
\qq\qq\qq\q\\
\ns\ds\le C\big(\|\d u\|_{L^{\infty}_\dbF(0,T;H_1)}\cdot\|\d u\|_{L^{8}_\dbF(0,T;H_1)}^2\big) ,
\end{array}
\end{equation}
\begin{equation}\label{th suff-eq13add}
\begin{array}{ll}\ds
\q \Big|
\dbE \int_{0}^{T}\Big(\big\langle\hat{g}_{22}(t)\d u(t), \d u(t)\big\rangle_{H_{1}}- \frac{1}{2}\big\langle g_{22}(t)\d u(t), \d u(t)\big\rangle_{H_{1}}\Big)dt\Big|
\qq\qq\qq\q\\
\ns\ds\le C\big(\|\d u\|_{L^{\infty}_\dbF(0,T;H_1)}\cdot\|\d u\|_{L^{8}_\dbF(0,T;H_1)}^2\big) ,
\end{array}
\end{equation}
and
\begin{equation}\label{th suff-eq14}
\begin{array}{ll}\ds
\q \Big|
\dbE \Big(\big\langle \hat{h}_{xx}(\bar x(T))\delta x(T),\delta x(T)\big\rangle_H -\frac{1}{2}\big\langle h_{xx}(\bar x(T))y(T), y(T)\big\rangle_H\Big)\Big|
\qq\qq\q \ \\
\ns\ds\le C\big(\|\d u\|_{L^{\infty}_\dbF(0,T;H_1)}\cdot\|\d u\|_{L^{8}_\dbF(0,T;H_1)}^2\big).
\end{array}
\end{equation}

Also, by  (\ref{th suff estmate r2}) and
\textbf{(A3)},
\begin{equation}\label{th suff-eq15}
\begin{array}{ll}\ds
\q \Big| \dbE \int_{0}^{T}\Big\langle g_{1}(t),
\d x(t)-y(t)-\frac{1}{2}z(t)\Big\rangle_Hdt
+\dbE~\Big\langle h_{x}(\bar x(T)),\d x(T)-y(T)-\frac{1}{2}z(T)\Big\rangle_H \Big| \\
\ns\ds\le \Big\|\d
x-y-\frac{1}{2}z\Big\|_{C_\dbF([0,T];L^{2}(\O;H))}
\le C\big(\|\d
u\|_{L^{\infty}_\dbF(0,T;H_1)}\cdot\|\d
u\|_{L^{8}_\dbF(0,T;H_1)}^2\big).
\end{array}
\end{equation}

Combining (\ref{th suff-eq10})--(\ref{th
suff-eq11})  with (\ref{th suff-eq12})--(\ref{th
suff-eq15}), we have that
\begin{equation}\label{th suff-eq16}
\begin{array}{ll}\ds
\q \cJ(u) - \cJ(\bar u)\\
\ns\ds \ge \mE\int_0^T \[\big\langle g_1(t),
y(t)\big\rangle_H + \frac{1}{2} \big\langle
g_1(t), z(t)\big\rangle_H + \big\langle
g_2(t), \d u(t)\big\rangle_{H_{1}} \\
\ns\ds \qq\qq + \frac{1}{2}\(\big\langle
g_{11}(t)y(t), y(t)\big\rangle_H +
2\big\langle g_{12}(t)y(t), \d
u(t)\big\rangle_{H_{1}} + \big\langle g_{22}(t)\d
u(t), \d u(t)\big\rangle_{H_{1}} \)
\]dt \\
\ns\ds\q + \mE\[ \big\langle h_x(\bar x(T)),
y(T) \big\rangle_H + \frac{1}{2}\big\langle
h_x(\bar x(T)), z(T) \big\rangle_H +
\frac{1}{2}\big\langle h_{xx}(\bar x(T))
y(T),y(T)\big\rangle_H
\]\\[+0.4em]
\ns\ds \q - C \|\d
u\|_{L^{\infty}_\dbF(0,T;H_1)}\cdot \|\d
u\|_{L^{8}_\dbF(0,T;H_1)}^2.
\end{array}
\end{equation}
Substituting
\eqref{th3.1-eq15}--\eqref{th3.1-eq17} into
\eqref{th suff-eq16} and combining with
\eqref{1st order condition}, we get that
\begin{equation}\label{Taylor expansion}
\begin{array}{ll}\ds
\q \cJ(u)-\cJ(\bar u)\\
\ns\ds \ge -\mE\int_0^T\big\langle\dbH_{u}(t),\d u(t)\big\rangle_{H_{1}}dt\\
\ns\ds\q
-\mE\int_0^T \(
\frac{1}{2}\big\langle \dbH_{uu}(t)\d
u(t), \d u(t) \big\rangle_{H_{1}}\! +\!
\frac{1}{2}\big\langle
b_2(t)^{*}P(t)b_2(t)\d u(t), \d
u(t)\big\rangle_{H_{1}}
\)dt\\
\ns\ds \q- \mE\!\int_0^T\!\!\!
\big\langle\big(\dbH_{xu}(t)\! +\! a_2(t)^*
P(t)\! +\! b_2(t)^*P(t) b_1(t)\big)y(t), \d
u(t)\big\rangle_{H_{1}} dt  \\
\ns\ds \q-\frac{1}{2}\mE\int_0^T
\big\langle \big(\widehat Q^{(0)}+ Q^{(0)}\big)(0,a_2(t)\d u(t),
b_2(t)\d u(t)), b_2(t)\d u(t)
\big\rangle_{H}dt\\
\ns\ds\q - C\big(\|\d u\|_{L^{\infty}_\dbF(0,T;H_1)}\cdot\|\d u\|_{L^{8}_\dbF(0,T;H_1)}^2\big)\\
\ns\ds= -\mE\int_0^T\big\langle\dbH_{u}(t),\d u(t)\big\rangle_{H_{1}}dt-\frac{1}{2}\Lambda(\d u)
- C\big(\|\d u\|_{L^{\infty}_\dbF(0,T;H_1)}\cdot\|\d u\|_{L^{8}_\dbF(0,T;H_1)}^2\big)\\
\ns\ds\ge \varrho\|\d u\|_{L^{8}_\dbF(0,T;H_1)}^2- C\big(\|\d u\|_{L^{\infty}_\dbF(0,T;H_1)}\cdot\|\d u\|_{L^{8}_\dbF(0,T;H_1)}^2\big).
\end{array}
\end{equation}
Then, choosing $\sigma$ small enough such that
$$C\|\d u\|_{L^{\infty}_\dbF(0,T;H_1)}\le C\sigma\le \frac{\varrho}{2},$$
we finally obtain (\ref{strong local min}).
\endpf

\vspace{0.2cm}

When the BSEE \eqref{op-bsystem3} has a unique
$V$-transposition solution $(P,Q)$ with which
$P_T$, $J(\cd)$, $K(\cd)$ and $F(\cd)$ are given
by (\ref{zv2}), the following result immediately
follows from Theorem \ref{th second order suff}.

\begin{corollary}\label{th 2nd suff v trans}
In addition to the assumptions in Theorem
\ref{th second order suff}, assume that
\textbf{(A5)} holds.  If there exists a constant
$\varrho>0$ such that for any $u(\cd) \in
\cU^{\infty}[0,T]$,
\begin{equation*}
\mE\int_{0}^{T}\big\langle\dbH_{u}(t), u(t)-\bar u(t)  \big\rangle_{H_{1}} dt\le0
\end{equation*}
and
\begin{equation*}
\begin{array}{ll}\ds
\tilde\Lambda(u-\bar u)
\leq - 2\varrho\|u-\bar u\|_{L^{8}_\dbF(0,T;H_1)}^2,
\end{array}
\end{equation*}
then there exists a constant $\sigma>0$ such that for any $u(\cd) \in \cU^{\infty}[0,T]$ with $\|u-\bar u\|_{L^{\infty}_\dbF(0,T;H_1)}\le \sigma$,
\begin{equation*}
\cJ(u)\ge \cJ(\bar u)+\frac{\varrho}{2}\|u-\bar u\|_{L^{8}_\dbF(0,T;H_1)}^2,
\end{equation*}
and, $\bar u$ is a local minima of the optimal control problem (\ref{jk1}).
\end{corollary}

In what follows we refine the second order
sufficient conditions in Theorem \ref{th second
order suff} and Corollary \ref{th 2nd suff v
trans} by the general Legendre form.

\begin{definition}
Let $X$ be a reflexive Banach space. A
functional $\Psi: X\to \dbR$ is called a general
Legendre form if $\Psi$ is weakly lower
semicontinuous, positively homogeneous of degree
2, i.e, for any $x\in X$, $t>0$,
$\Psi(tx)=t^2\Psi (x)$ and if
$x_{k}\weakconvergent{w}x$  and $\Psi(x_{k})\to
\Psi(x)$, it holds that $x_{k}\to x$ strongly.
\end{definition}

Some sufficient and necessary conditions to
ensure a functional to be a  Legendre form can
be found in \cite{Bonnans00}.

Define
$$
T_{\cU^{8}[0,T]}(\bar u)\=cl_{8}
\Big\{v=\alpha (u-\bar u)\ |\  u\in
\cU^{\infty}[0,T], \alpha\ge 0\Big\},
$$
where $cl_{8}(\cA)$ is the closure of a set
$\cA$ under the norm topology of the Banach
space $L^{8}_\dbF(0,T;H_1)$. If $-\Lambda$ is a
general Legendre form defined on
$L^{8}_\dbF(0,T;H_1)$, the negative definite
condition  (\ref{2nd order condition}) can be
weaken to the following directional negative
definite condition:

\no{\bf (A9)} \vspace{-0.4cm}
\begin{equation}\label{directional NDF}
\Lambda(v) <0,\q \forall v\in \mathcal{C}_{ \cU^{8}[0,T]}(\bar u),
\end{equation}
where
$$
\mathcal{C}_{\cU^{8}[0,T]}(\bar u)\=\[
T_{\cU^{8}[0,T]}(\bar
u)\bigcap\Big\{v\in\cU^{8}[0,T]\; \Big|\
\mE\int^{T}_{0}\!\!\inner{\dbH_{u}(t)}{v(t)}_{H_{1}}dt=0
\Big\}\]\setminus \{0\}.
$$
When  $-\Lambda$ is a general Legendre form, we have the following second order sufficient condition:

\begin{theorem}\label{th second order suff general}
Assume that $x_0 \in  L^2_{\cF_0}(\O;H)$ and
$-\Lambda$ is a general Legendre form  on
$L^{8}_\dbF(0,T;$ $H_1)$. Let
{\bf(A1)}--{\bf(A3)} and {\bf(A7)}--{\bf(A9)}
hold, and let $\bar u(\cd)$ be an admissible
control and $\bar x(\cd)$ be the corresponding
state. If for any $u(\cd) \in
\cU^{\infty}[0,T]$,
\begin{equation}\label{12.12-eq2}
\mE\int_{0}^{T}\big\langle\dbH_{u}(t), u(t)-\bar
u(t)  \big\rangle_{H_{1}} dt\le0,
\end{equation}
then there exist constants $\sigma>0$ and
$\varrho>0$ such that for any $u(\cd) \in
\cU^{\infty}[0,T]$ with $\|u-\bar
u\|_{L^{\infty}_\dbF(0,T;H_1)}\le \sigma$, the
quadratic growth condition (\ref{strong local
min}) holds.
\end{theorem}

{\it Proof}\,: We prove this conclusion through
a contradiction argument. If one could not find
$\sigma>0$ and $\varrho>0$ such that
(\ref{strong local min}) holds, then there must
exist sequences $\{\varrho_{n}\}_{n}^{\infty}$
and $\{u_{n}\}_{n=1}^{\infty}$ such that for any
$n$, $\varrho_{n}>0$, $\varrho_{n}\to 0$,
$u_{n}\in \cU^{\infty}[0,T]$, $\|u_{n}-\bar
u\|_{L^{\infty}_\dbF(0,T;H_1)}\to 0$ (as $n\to
\infty$) and
\begin{equation}\label{contr quadratic growth}
\cJ(u_{n})< \cJ(\bar u)+\frac{\varrho_{n}}{2}\|u_{n}-\bar u\|_{L^{8}_\dbF(0,T;H_1)}^2.
\end{equation}
Let
$$
v_{n}=\frac{u_{n}-\bar u}{\|u_{n}-\bar u\|_{L^{8}_\dbF(0,T;H_1)}}.
$$
It is clear that $v_{n}$ is a unit vector of
$L^{8}_\dbF(0,T;H_1)$ for any $n\in\dbN$,  and
there exists a subsequence
$\{u_{n_{k}}\}_{k=1}^{\infty}$ which converges
weakly to a vector $v\in L^{8}_\dbF(0,T;H_1)$.
Without loss of generality, we assume
$v_{n}\weakconvergent{w} v$. Since $U$ is
convex, $T_{\cU^{8}[0,T]}(\bar u)$ is a closed
convex set. Noting that $L^{8}_\dbF(0,T;H_1)$ is
a reflexive Banach space, we have $v\in
T_{\cU^{8}[0,T]}(\bar u)$.

\vspace{0.2cm}

Let us divide the rest of the proof into three
steps.

\vspace{0.2cm}

{\bf Step 1:} In this step, we prove that
$\ds\mE\int_{0}^{T}\inner{\dbH_{u}(t)}{v(t)}dt=0$.
By \eqref{12.12-eq2},
$$
\mE\int_{0}^{T}\inner{\dbH_{u}(t)}{v_{n}(t)}_{H_{1}}dt\le 0,
$$
and hence
$$\mE\int_{0}^{T}\inner{\dbH_{u}(t)}{v(t)}_{H_{1}}dt\le 0.$$
If for some $\varepsilon>0$,
$$\mE\int_{0}^{T}\inner{\dbH_{u}(t)}{v}_{H_{1}}dt<-\varepsilon<
0,
$$
by \eqref{Taylor expansion}, it is easy to find
that
\begin{equation*}
\3n\3n\begin{array}{ll}\ds \q \cJ(u_{n})
\3n&\ds\ge \cJ(\bar
u)-\mE\int_0^T\big\langle\dbH_{u}(t),u_{n}(t)-\bar
u(t)\big\rangle_{H_{1}}dt
- o(\|u_{n}(t)-\bar u(t)\|_{L^{8}_\dbF(0,T;H_1)})\\
\ns&\ds = \cJ(\bar u)\!-\!\|u_{n}\!-\!\bar
u\|_{L^{8}_\dbF(0,T;H_1)}
\mE\int_0^T\big\langle\dbH_{u}(t),v_{n}(t)\big\rangle_{H_{1}}dt
\!-\! o(\|u_{n}(t)\!-\!\bar u(t)\|_{L^{8}_\dbF(0,T;H_1)}).
\end{array}
\end{equation*}
Then, by assumption (\ref{contr quadratic
growth}), we have that
$$
-\mE\int_0^T\big\langle\dbH_{u}(t),v_{n}(t)\big\rangle_{H_{1}}dt
- o(1)< \frac{\varrho_{n}}{2} \|u_{n}-\bar u\|_{L^{8}_\dbF(0,T;H_1)}.
$$
Letting $n\to \infty$, we get that
$$0<\varepsilon< -\mE\int_0^T\big\langle\dbH_{u}(t),v(t)\big\rangle_{H_{1}}dt
\le 0,$$ a contradiction. Therefore,
$$\mE\int_{0}^{T}\inner{\dbH_{u}(t)}{v(t)}dt=0.$$

\ms

{\bf Step 2:} In this step, we prove that $v\neq
0$.

\ms

If not, $\Lambda (v)=\Lambda(0)=0$. Using
(\ref{Taylor expansion}) again, we obtain that
\begin{equation*}
\!\!\!\! \begin{array}{ll}\ds \cJ(u_{n})
\3n&\ds\ge \cJ(\bar
u)-\mE\int_0^T\big\langle\dbH_{u}(t),u_{n}(t)-\bar
u(t)\big\rangle_{H_{1}}dt
-\frac{1}{2}\Lambda(u_{n}(t)-\bar u(t))\\
\ns&\ds\q
- C\big(\|u_{n}(t)-\bar u(t)\|_{L^{\infty}_\dbF(0,T;H_1)} \|u_{n}(t)-\bar u(t)\|_{L^{8}_\dbF(0,T;H_1)}^2\big)\\
\ns&\ds \ge \cJ(\bar
u)-\frac{1}{2}\Lambda(u_{n}(t)-\bar u(t))
- C\big(\|u_{n}(t)-\bar u(t)\|_{L^{\infty}_\dbF(0,T;H_1)} \|u_{n}(t)-\bar u(t)\|_{L^{8}_\dbF(0,T;H_1)}^2\big).\\
\end{array}
\end{equation*}
By \eqref{contr quadratic growth},
$$
\!\!\begin{array}{ll}\ds
-\!\frac{1}{2}\Lambda(u_{n}(t)\!-\!\bar u(t))
\!-\! C\big(\|u_{n}(t)\!-\!\bar
u(t)\|_{L^{\infty}_\dbF(0,T;H_1)}
\|u_{n}(t)\!-\!\bar
u(t)\|_{L^{8}_\dbF(0,T;H_1)}^2\big)\! <\!
\frac{\varrho_{n}}{2} \|u_{n}\!-\!\bar
u\|_{L^{8}_\dbF(0,T;H_1)}^2,
\end{array}
$$
which implies
\begin{equation}\label{Th4.2 eq1}
-\frac{1}{2}\Lambda(v_{n})
- C\|u_{n}(t)-\bar u(t)\|_{L^{\infty}_\dbF(0,T;H_1)}< \frac{\varrho_{n}}{2}.
\end{equation}
Since $-\Lambda$ is weakly lower semicontinuous,
$$ \liminf_{n\to \infty}-\Lambda(v_{n})\ge -\Lambda(v).$$
Then, by (\ref{Th4.2 eq1}),
$$
0=-\frac{1}{2}\Lambda(v)\le  \liminf_{n\to \infty}\(-\frac{1}{2}\Lambda(v_{n}) - C\|u_{n}(t)-\bar u(t)\|_{L^{\infty}_\dbF(0,T;H_1)}\)\le\liminf_{n\to \infty} \frac{\varrho_{n}}{2}=0,
$$
which implies that there exists a subsequence $\{v_{n_{k}}\}_{k=1}^{\infty}$ (of $\{v_{n}\}_{n=1}^{\infty}$) such that
$-\Lambda (v_{n_{k}})\to -\Lambda (v)=0$ as $k\to \infty$. Since $-\Lambda$ is a Legendre form and $v_{n_{k}}\weakconvergent{w}v$, we have $v_{n_{k}}\to v$ strongly. But $\|v_{n_{k}}\|_{L^{8}_\dbF(0,T;H_1)}=1$ and therefore $\|v\|_{L^{8}_\dbF(0,T;H_1)}=1$, contradicting to the assumption that $v=0$.

\vspace{0.2cm}

{\bf Step 3:} By Steps 1 and 2, we have proved
that $v\in \mathcal{C}_{ \cU^{8}[0,T]}(\bar u)$.
Then by Assumption \textbf{(A8)}, there exists a
constant $\varepsilon>0$ such that
$$
-\Lambda (v)\ge \varepsilon>0,
$$
which gives
$$0<\frac{\varepsilon}{2}\le -\frac{1}{2}\Lambda (v)\le \liminf_{n\to \infty}\(-\frac{1}{2}\Lambda(v_{n}) - C\|u_{n}(t)-\bar u(t)\|_{L^{\infty}_\dbF(0,T;H_1)}\)\le\liminf_{n\to \infty} \frac{\varrho_{n}}{2}=0,$$
a contradiction. This completes the proof of
Theorem \ref{th second order suff general}.
\endpf

\begin{corollary}
Assume that $x_0\in L^2_{\cF_0}(\O;H)$,
{\bf(A1)}--{\bf(A3)},  {\bf(A5)} and
{\bf(A7)}--{\bf(A8)} hold. Let $(P,Q)$ be the
unique $V$-transposition solution to BSEE
\eqref{op-bsystem3} with $P_T$, $J(\cd)$,
$K(\cd)$ and $F(\cd)$ given by (\ref{zv2}) and
let $\bar u(\cd)$ be an admissible control with
$\bar x(\cd)$ the corresponding state. If
$-\tilde \Lambda$ is a Legendre form on
$L^{8}_\dbF(0,T;H_1)$,
$$
\tilde \Lambda(v)< 0,\q \forall\ v\in \mathcal{C}_{ \cU^{8}[0,T]}(\bar u)
$$
and for any $u(\cd) \in \cU^{\infty}[0,T]$,
\begin{equation*}
\mE\int_{0}^{T}\big\langle\dbH_{u}(t), u(t)-\bar u(t)  \big\rangle_{H_{1}} dt\le0,
\end{equation*}
then there exist  constants $\sigma>0$ and
$\varrho>0$ such that for any $u(\cd) \in
\cU^{\infty}[0,T]$ with $\|u-\bar
u\|_{L^{\infty}_\dbF(0,T;H_1)}\le \sigma$, the
quadratic growth condition (\ref{strong local
min}) holds.
\end{corollary}

\begin{remark}
The proof of Theorem \ref{th second order suff
general} is a modification  of the related
conclusion in deterministic optimization
problem, see \cite[Chapter 3]{Bonnans00}.  The
corresponding results in the deterministic
optimal control problem can be found in
\cite{BDP2}, and that in optimal control
problems of stochastic differential equations
can be found in \cite{BS}. Note that, in Theorem
\ref{th second order suff general}, we do not
need the assumptions that $ b_{uu}(t,x,u)\equiv
0$  or the maps $(x,u)\mapsto a(t,x,u)$ and
$(x,u)\mapsto b(t,x,u)$ are affine for a.e.
$t\in [0,T]$. Therefore, Theorem \ref{th second
order suff general} is much more general than
\cite[Proposition 4.15]{BS}. In addition, even
though the condition (\ref{directional NDF}) in
Theorem \ref{th second order suff general} is
weaker than the condition (\ref{2nd order
condition}) in Theorem \ref{th second order
suff}, in the stochastic cases, there exist some
essential difficulties to verify if the
corresponding quadratic-like forms $-\Lambda$ or
$-\tilde\Lambda$ are Legendre form (see \cite[ Examples
4.16--4.17]{BS}). Therefore, sometimes
it is much more convenient to use Theorem
\ref{th second order suff} in practice.
\end{remark}

\section{Examples}\label{s6}

In this section, we shall give some examples.
Firstly, we apply our second order necessary
condition for systems of controlled stochastic
heat equations. The same thing can be done for
lots of other systems, such as stochastic
Schr\"odinger equations,  stochastic Korteweg-de
Vries equations, stochastic Kuramoto-Sivashinsky
equations, stochastic Cahn-Hilliard equations,
etc.

\begin{example}
Let $H=L^{2}[0,1]\times L^{2}[0,1]$,
$H_1=H_0^1(0,1)\times H_0^1(0,1)$,
$V=H^{-1}(0,1)\times H^{-1}(0,1)$ and
$U=H_0^1(0,1)\times B_{H_0^1(0,1)}$ where $
B_{H_0^1(0,1)}$ is the closed unit ball in
$H_0^1(0,1)$. Then $V'=H_0^1(0,1)\times
H_0^1(0,1)$. Define an operator $A$ by
$$
\begin{cases}\ds
D(A)=H^2(0,1)\cap H_0^1(0,1),\\[-0.2em]
\ns\ds Af=\pa_{xx}f,\q\forall f\in D(A).
\end{cases}
$$
It is clear that the embedding from $H$ to $V$
is Hilbert-Schmidt and $A$ generates a
$C_0$-semigroup on $H_0^1(0,1)$.  Consider the
following control system:
\begin{equation}\label{example control system}
\left\{
\begin{array}{lll}
\ds d\f_{1} = \pa_{xx}\f_{1} dt +  u_{1}dt +(\f_{1}+\f_{2}) dW(t) &\mbox{ in }(0,T]\times (0,1),\\[+0.4em]
\ds d\f_{2} =\pa_{xx}\f_{2} dt+ u_{2}^2 dW(t) &\mbox{ in }(0,T]\times (0,1),\\[+0.4em]
\ns\ds \f_{1}(t,0)=\f_{1}(t,1)=0, &\mbox{ in }
(0,T],\\
\ns\ds \f_{2}(t,0)=\f_{2}(t,1)=0, &\mbox{ in }
(0,T],\\
\ns\ds \f_{1}(0,x) =\phi(x), &\mbox{ on } (0,1)\\
\ns\ds \f_{2}(0,x) =0, &\mbox{ on } (0,1),
\end{array}
\right.
\end{equation}
and the cost functional
$$
\cJ(u)=\frac{1}{2}\mE \langle
\f_{1}(T),\f_{1}(T)\rangle_{L^2(0,1)}.
$$

It is easy to see that {\bf(A1)}--{\bf(A3)} hold
for the above optimal control problem.
Furthermore,
$$
a_x=0\in L^\infty_\dbF(0,T;\cL_{HV'}),\;
b_x=\left(
\begin{array}{cc}
I & I \\
0 & 0 \\
\end{array}
\right)\in L^\infty_\dbF(0,T;\cL_{HV'}).
$$
Then, we see that {\bf(A5)} holds.

Let $\ds\phi=\sum_{n=1}^\infty a_n \sqrt{2}\sin
n\pi x\in L^2(0,1)$. We claim that
$(u_1,u_2)=(f,0)\in \mathcal{U}^2[0,T]$, where
$$
f(t,x)=\sum_{n=1}^\infty f_n(t)\sqrt{2}\sin n\pi
x\q \mbox{ for } f_n(t) =
-\frac{a_n}{T}e^{-(n^2\pi^2+1/2)t + W(t)}
$$
is an optimal control. Indeed, if
$(u_1,u_2)=(f,0)$, then the corresponding
solution $(\f_1,\f_2)$ satisfies that
$\f_{1}(T)=0$. Next, direct computations show
that $f\in L^2_{\dbF}(0,T;H^1_0(0,1))$. This
verifies our claim. Furthermore, one can show
that
$f(\cdot)\in\dbL_{2,\dbF}^{1,2}(H^1_0(0,1))$.
Hence, we find that the first condition in
{\bf(A6)} holds.

For this optimal control problem, the
Hamiltonian is
$$
\dbH(t,(\f_{1},\f_{2}),(u_{1},u_{2}),(p_{1},p_{2}),(q_{1},q_{2}))=p_{1}u_{1}+
q_{1}(\f_{1}+\f_{2})+ q_2u_{2}^2,
$$
and the
corresponding first order adjoint
equation is
\begin{equation}\label{example 1st adj system}
\left\{
\begin{array}{lll}
\ds dp_{1} = -\pa_{xx}p_{1} dt-q_{1}dt+ q_{1}dW(t) &\mbox{ in }[0,T)\times(0,1),\\[+0.4em]
\ds dp_{2} = -\pa_{xx}p_{2} dt-q_{1}dt+ q_{2}dW(t) &\mbox{ in }[0,T)\times(0,1),\\
\ns\ds p_{1}(\cd,0)=p_{1}(\cd,1)=0, &\mbox{ on } [0,T),\\
\ns\ds p_{2}(\cd,0)=p_{2}(\cd,1)=0, &\mbox{ on } [0,T),\\
\ns\ds p_{1}(T,\cd) =p_{2}(T,\cd)=0, &\mbox{ in
}(0,1).
\end{array}
\right.
\end{equation}
Obviously, $(p_{1},p_{2})\equiv 0$,
$(q_{1},q_{2})\equiv 0$ and therefore,
$$
\dbH(t,(\f_{1},\f_{2}),(u_{1},u_{2}),(p_{1},q_{1}),(p_{2},q_{2}))\equiv
0.
$$
Then, the second order adjoint equation reads
\begin{equation}\label{example 2nd adj system}
\left\{\!\!\!
\begin{array}{lll}
\ds dP = -\left(
          \begin{array}{cc}
          A^{*}P_{11}+P_{11}A^{*} & A^{*}P_{12}+P_{12}A^{*} \\
          A^{*}P_{21}+P_{21}A^{*} & A^{*}P_{22}+P_{22}A^{*} \\
          \end{array}
          \right) dt \\[+1em]
\;\;\;\qq - \bigg[\left(
\begin{array}{cc}
P_{11} & P_{11} \\
P_{11} & P_{11} \\
\end{array}
\right)
\!+\! \left(\!
\begin{array}{cc}
2Q_{11} & Q_{11}\!+\!Q_{12} \\
Q_{11}\!+\!Q_{21} & Q_{12}+\!Q_{21} \\
\end{array}
\right)\bigg]dt\! +\! \left(
\begin{array}{cc}
Q_{11} & Q_{12} \\
Q_{21} & Q_{22} \\
\end{array}
\!\right)dW(t)\;\, \mbox{in }[0,T),\\[+1em]
\ns\ds P(T)=\left(
\begin{array}{cc}
-I & 0\\
0 & 0 \\
\end{array}
\right).
\end{array}
\right.
\end{equation}
It is clear that $Q=\left(
\begin{array}{cc}
Q_{11} & Q_{12} \\
Q_{21} & Q_{22} \\
\end{array}
\right)=0$. Then $P$ is the solution to
$$
\left\{
\begin{array}{lll}
\ds dP = -\left(
          \begin{array}{cc}
          A^{*}P_{11}+P_{11}A^{*} & A^{*}P_{12}+P_{12}A^{*} \\
          A^{*}P_{21}+P_{21}A^{*} & A^{*}P_{22}+P_{22}A^{*} \\
          \end{array}
          \right) dt
- \left(
\begin{array}{cc}
P_{11} & P_{11} \\
P_{11} & P_{11} \\
\end{array}
\right)dt&\mbox{ in }[0,T),\\[+1em]
\ns\ds P(T)=\left(
\begin{array}{cc}
-I & 0\\
0 & 0 \\
\end{array}
\right).
\end{array}
\right.
$$
Obviously,
$$
\begin{array}{ll}\ds
P(\cd) =-e^{A(T-\cd)}Je^{A(T-\cd)} - \int_\cd^T
e^{A(s-\cd)}P(s)e^{A(s-\cd)}ds\\
\ns\ds\in
\dbL_{2,\dbF}^{1,2}\big(\cL_2[(H^1_0(0,1))^2;(L^2(0,1))^2]\big)\cap
L^{\infty}\big([0,T]\times\Omega;\cL_2[(H^1_0(0,1))^2;(L^2(0,1))^2]\big).
\end{array}
$$
Further, by the classical theory of Riccati
equations (see \cite[Part IV, Section 2.2,
Theorem 2.1]{Bensoussan3}), we know that
$P_{11}(t)<0$, for any $t\in [0,T]$.

Since $b_u(t,\bar x(t),\bar u(t))\equiv 0$, we have,
$$
\dbH_{u}(t,(\f^{1},\f^{2}),(u_{1},u_{2}),(p^{1},q^{1}),(p^{2},q^{2}))\equiv0,
$$
$$
\dbH_{uu}(t,(\f^{1},\f^{2}),(u_{1},u_{2}),(p^{1},q^{1}),(p^{2},q^{2}))\! +\!
b_u(t,\bar x(t),\bar u(t))^{*}P(t)b_u(t,\bar
x(t),\bar u(t))\equiv0,
$$
and
$$
\dbS(t)=\left(
          \begin{array}{cc}
            P_{11} & P_{12} \\
            0 & 0 \\
          \end{array}
        \right), \q\nabla \dbS(t)=0.
$$
Then,  we see that the second condition in
{\bf(A6)} holds.

\end{example}

In what follows, we consider an application of the second order sufficient condition in the stochastic LQ problems.
\begin{example}

Let us consider the following linear control system
\begin{eqnarray}\label{Linear system}
\left\{
\begin{array}{lll}\ds
dx = \big(Ax +B_{1}x+C_{1}u\big)dt +\big(B_{2}x+C_{2}u \big) dW(t) &\mbox{ in }(0,T],\\
\ns\ds x(0)=x_0,
\end{array}
\right.
\end{eqnarray}
and the cost functional
\begin{equation}\label{LQ cost Lagrange}
J(u(\cdot))\!=\!\frac{1}{2}\dbE\int_{0}^{T}\!\!\(\inner{Rx(t)}{x(t)}_H
+2\inner{x(t)}{Mu(t)}_H+\inner{Nu(t)}{u(t)}_{H_{1}}\!\)dt.
\end{equation}
We assume that
$B_{1}, B_{2}\in \cL_{HV^{'}}$, $C_{1},C_{2},M\in \cL(H_1;H)$, $R\in \cL(H)$, $N\in \cL(H_1)$. Moreover, $R$ and $N$  are self-adjoint.

It is clearly that, the optimal control problem (\ref{jk2}) for control system (\ref{Linear system}) and cost functional (\ref{LQ cost Lagrange}) is well-defined on $\cU^{2}[0,T]$.

Let $(\bar x, \bar u)$ be an admissible  pair. Define
\begin{equation*}
\begin{array}{ll}\ds
\dbH(t,x,u,k_1,k_2) \= \big\langle k_1, B_{1}x+C_{1}u  \big\rangle_H + \big\langle k_2, B_{2}x+C_{2}u  \big\rangle_H\\
\ns\ds\hspace{3cm}\q
-\frac{1}{2}\big(\big\langle  R x(t),x(t)\big\rangle_H
+2\big\langle x(t),  M u(t)\big\rangle_H+\big\langle  N u(t), u(t)\big\rangle_{H_{1}}\! \big) ,\\
\ns\ds \hspace{5.16cm} (t,x,u,k_1,k_2)\in
[0,T]\times H \times U\times H\times H,
\end{array}
\end{equation*}
and define the first and second order adjoint equations:
\begin{equation}\label{1st adjoint LQ}
\left\{
\begin{array}{ll}
dp=-A^{*}pdt -\big(B_{1}^{*}p+  B_{2}^{*}q-
R\bar{x}(t)-  M\bar{u}(t)\big)dt+qdW(t) &\mbox{
in }[0,T),
\\
p(T)=0
\end{array}
\right.
\end{equation}
and
\begin{equation}\label{2nd adjoint LQ}
\left\{
\begin{array}{ll}
dP=-\big(A^{*}+B^{*}_{1}\big)Pdt-P\big(A+B_{1}\big)dt - B_{2}^{*}P B_{2}dt\\[+0.4em]
\qq\q\;\; -\big(B_{2}^{*}Q+ QB_{2}\big)dt+ R dt+QdW(t) &\mbox{ in } [0,T), \qq\\
P(T)=0.
\end{array}\right.
\end{equation}
Obviously, BSEE \eqref{1st adjoint LQ} admits a
unique transposition solution $\big(p, q\big)$,
and, BSEE \eqref{2nd adjoint LQ} admits a
unique $V$-transposition solution $\big(P,
Q\big)$. In addition, since the operators
$B_{1},B_{2}$ and $N$ independent of
$(t,\omega)$, we have that $\big(P, Q\big)$ is
actually the solution to the follow
deterministic operator-valued evolution
equation:
\begin{equation}\label{2nd adjoint LQ deterministic}
\left\{
\begin{array}{l}
dP=-\big(A^{*}+B^{*}_{1}\big)Pdt-P\big(A+B_{1}\big)dt - B_{2}^{*}P B_{2}dt+ Rdt\qq \mbox{in}\; [0,T), \\
P(T)=0.
\end{array}\right.
\end{equation}

Let $u(\cdot)\in \cU^{2}[0,T]$ be another admissible control with the corresponding state $x(\cdot)$ and denote $\delta x= x(\cdot)-\bar x(\cdot)$, $\delta u= u(\cdot)-\bar u(\cdot)$. We have that
\begin{equation}\label{Taylor expansion for LQ}
\begin{array}{ll}\ds
\q \cJ(u)-\cJ(\bar u)\\
\ns\ds =\mE\int_0^T\( \big\langle R \bar x(t), \delta x(t) \big\rangle_{H}
+\frac{1}{2}\big\langle R \delta x(t), \delta x(t) \big\rangle_{H}\\
\ns\ds\q
+\big\langle N \bar u(t), \delta u(t) \big\rangle_{H_{1}}
+\frac{1}{2}\big\langle N \delta u(t), \delta u(t) \big\rangle_{H_{1}}\\
\ns\ds\q
+\big\langle \bar x(t),   M \delta u(t)  \big\rangle_{H}
+\big\langle \delta  x(t),  M \bar u(t)  \big\rangle_{H}
+\big\langle \delta  x(t),   M \delta u(t) \big\rangle_{H}\)dt
\\
\ns\ds
=-\mE\int_0^T \(
\big\langle \dbH_{u}(t), \d u(t) \big\rangle_{H_{1}}
-\frac{1}{2}\big\langle R \delta x(t), \delta x(t) \big\rangle_{H}\\
\ns\ds
\q-\frac{1}{2}\big\langle  N \delta u(t), \delta
u(t) \big\rangle_{H_{1}} -\big\langle \delta
x(t),   M \delta u(t) \big\rangle_{H}\)dt.
\end{array}
\end{equation}
By It\^{o}'s formula,
\begin{eqnarray}\label{equ LQ and 2nd condition}
&&\q-\mE\int_0^T \( \big\langle R \delta x(t),
\delta x(t) \big\rangle_{H}+\big\langle  N
\delta u(t), \delta u(t) \big\rangle_{H_{1}}
+2\big\langle \delta  x(t),   M \delta u(t) \big\rangle_{H}\)dt\nonumber\\[+0.5em]
&& =-\mE\int_0^T\[ 2\big\langle \dbS \delta
x(t), \delta u(t)
\big\rangle_{H_{1}}+\big\langle \big(N
+C_{2}^{*}P_{2}C_{2}\big) \delta u(t), \delta
u(t) \big\rangle_{H_{1}}
\]dt\nonumber\\[+0.5em]
&&=-\tilde \Lambda(\delta u(\cdot)),
\end{eqnarray}
where
$$
\dbS= C_{1}^{*}P+ C_{2}^{*}PB_{2}- M^{*}.
$$
Noting that in the this special case, the quadratic form $\tilde\Lambda$ can be extended into the Hilbert space ${L^{2}_\dbF(0,T;H_1)}$.
Therefore, the second order sufficient condition
\begin{equation}\label{quadratic form estimats LQ}
\tilde \Lambda(\delta u(\cdot))\le - \varrho \|\d u\|_{L^{2}_\dbF(0,T;H_1)}^2
\end{equation}
holds true if and only if the quadratic
functional (defined on ${L^{2}_\dbF(0,T;H_1)}$)
\begin{equation}\label{equ 2nd condition}
F(\delta u)\=-\dbE\int_{0}^{T}
\inner{\big(\Gamma^{*}R\Gamma+\Gamma^{*}M+N
\big)\delta u(t) }{\delta u(t)}_{H_{1}}dt\le -
\varrho \|\d u\|_{L^{2}_\dbF(0,T;H_1)}^2,
\end{equation}
where $\Gamma\delta u=\delta x$.

Using a similar argumentation as in Section \ref{s5}, we have, when condition (\ref{equ 2nd condition}) is satisfied,  any $(\bar x,\bar u)$ is an local optimal pair if
$$
\mE\int_0^T \( \big\langle\dbH_{u}(t), u(t)-\bar
u(t) \big\rangle_{H_{1}}dt\le 0,\q \forall\;
u(\cdot)\in \cU^{2}[0,T].
$$
Furthermore, since the inequality
\eqref{quadratic form estimats LQ} holds true
for any $u(\cdot)\in \cU^{2}[0,T]$, we know that
the  $(\bar x,\bar u)$ satisfying the above
inequality is the unique globally minimizer.
\end{example}
%

%\section*{Acknowledgement}

\end{document}